%% file: main.tex
\documentclass[oneside]{article}
\usepackage[a4paper,
    left=1in,      %
    right=1in,     %
    top=1.2in,     %
    bottom=1.2in   %
]{geometry}
\usepackage{bm}
\usepackage{amsmath}
\usepackage{hyperref}
\usepackage{multicol}
\usepackage{amsfonts}
\usepackage{graphicx}
\usepackage{amsthm}
\usepackage{cleveref}
\usepackage{subcaption}
\usepackage{mathrsfs}
\usepackage{authblk}
\usepackage{float}
\usepackage{xcolor}
\usepackage{amssymb} %
\usepackage{amsmath} %
\usepackage{xcolor}
\usepackage{bm} %
\usepackage{newtxtext}
\usepackage{newtxmath}
\usepackage{graphicx}
\usepackage{subcaption}
\usepackage{sidecap} %
\usepackage{soul}

\usepackage[LGR,T1]{fontenc}
\usepackage{amsmath, amsthm, amssymb, verbatim, comment, tikz, pgfplots, mathtools, mathrsfs, ifthen, import, float, subcaption, enumitem, todonotes, bbm}
\tikzset{>=latex} %
\usepackage{tikz-network}
\usetikzlibrary{decorations.pathreplacing}
\usepackage{pgfplots}

\usepackage{amsfonts}
\usepackage[ruled,vlined,linesnumbered]{algorithm2e}%

\usepackage{parskip} %

\usepackage{xifthen}
\usepackage{enumerate}
\usepackage{titlesec}
\newcounter{reviewer}
\newcounter{point}[reviewer]
\newtheorem{theorem}{Theorem}[section]
\newtheorem{remark}[theorem]{Remark}
\newtheorem{proposition}[theorem]{Proposition}

\newtheorem{corollary}[theorem]{Corollary}

\makeatletter
\newcommand{\Romanum}[1]{\expandafter\@slowromancap\romannumeral #1@} %
\makeatother
\usepackage{graphicx} %
\usepackage{tikz}
\usetikzlibrary{patterns}
\tikzset{>=latex} %
\usepackage{tikz-network}
\usetikzlibrary{decorations.pathreplacing}
\usepackage{subcaption}
\title{
\textsc{Exact Integration of ReLU-type Neural Networks for Solving PDEs}}

\author{Shima Baharlouei\textsuperscript{1}\thanks{Corresponding author, email: shima.baharlouei@ehu.eus}\hspace{0.2cm},~Ignacio Brevis\textsuperscript{2},~Reuben Menezes\textsuperscript{2},~~David Pardo\textsuperscript{1, 3, 4},\\ \vspace{-2mm}~Kristoffer G. van der Zee\textsuperscript{2}\\%
\vspace{5mm}
\textsuperscript{1}{\small{University of the Basque Country (UPV/EHU), Leioa, Spain.}}\\%
\textsuperscript{2}{\small{Department of Mathematics, University of Nottingham, Nottingham, United Kingdom.}}\\%
\textsuperscript{3}{\small{Basque Center for Applied Mathematics (BCAM),  Bilbao, Spain.}}\\%
\textsuperscript{4}{\small{Ikerbasque (Basque Foundation for Sciences), Bilbao, Spain.}}}%

\date{\today}
\input{StandardSetOfCommands.tex}

\begin{document}
\maketitle
{\renewcommand{\thefootnote}{}\footnotetext{Emails of coauthors: $\{$ignacio.brevis1,kg.vanderzee, Reuben.Menezes$\}$@nottingham.ac.uk, david.pardo@ehu.eus}}
\begin{abstract}
Developing robust training methodologies for solving partial differential equations (PDEs) with neural networks is crucial for reliable scientific computing. Piecewise-linear activations such as the Rectified Linear Unit (ReLU) are a convenient choice in this setting as their structure resembles classical piecewise-linear approximations. 
However, when the loss functional incorporates the spatial gradient of the network, as in the Deep Ritz method, numerical integration becomes problematic. Specifically, relying on standard integration techniques (e.g., Monte Carlo sampling) during training yields incorrect gradients, even as the number of integration points approaches infinity, leading to optimization failure.
\par
\noindent To overcome this limitation, we present an  integration framework for single-hidden-layer ReLU networks that computes correct gradients when using Deep Ritz for PDEs or PINNs for first-order problems.
Our approach employs the Sutherland–Hodgman algorithm to partition the physical domain into subdomains according to breaking hyperplanes underlying the ReLU network. 
By finding these hyperplanes exactly, our method can correctly integrate the gradients, ensuring the optimizer receives correct information. %
\noindent In this work, we further assemble the energy in matrix--vector form that allows all neuron interactions to be computed efficiently with vectorized operations, introduce an alternative approach termed Conforming Monte Carlo that is a stochastic version of the proposed method, and propose two stochastic training strategies to reduce complexity. Numerical results on two-dimensional domains demonstrate that our framework recovers optimal PDE solutions that standard stochastic sampling and other existing quadrature methods fail to capture. 
\par
\noindent While the proposed approach is demonstrated here for elliptic PDEs within the Deep Ritz method using a single-hidden-layer network, the framework can apply broadly to variational PDE formulations involving ReLU networks.

\end{abstract}
\pagestyle{myheadings}
\thispagestyle{plain}
\markboth{\small Integration of ReLU for PDEs}{\small Integration of ReLU for PDEs}

\section{Introduction}

In recent years, neural networks (NNs) have emerged as a transformative computational tool for solving Partial Differential Equations (PDEs) \cite{baharlouei2026least, baharlouei2025least, raissi2019physics}. In NNs, researchers frequently employ the Rectified Linear Unit (ReLU) activation function due to its computational efficiency and effectiveness in mitigating vanishing gradients \cite{giovannelli2026limitation, goodfellow2016deep, he2015delving, he2020relu, opschoor2020deep}. Moreover, because the ReLU structure closely resembles classical piecewise-linear approximations, it is a convenient choice in this setting, which is attractive because of its familiarity with classical FEM. However, applying ReLU networks to solve PDEs via methods like the Deep Ritz Method~\cite{yu2018deep}, or PINNs for first-order problems~\cite{rowan2025failure}, introduces a challenge in numerical integration. Indeed, the choice of quadrature rule affects the accuracy of NN PDE solvers~\cite{rivera2022quadrature, taylor2026stochastic}. When using standard numerical integration (e.g., Monte Carlo) on the involved loss function, the computed gradients will converge towards an incorrect solution even as we increase the number of integration points, which may cause the optimizer to fail to converge to a correct minimizer.

\par
To illustrate this problem, consider the Deep Ritz Method for which the integrand of the loss contains spatial derivatives of the PDE solution. 
Computing the network's spatial gradients to evaluate the energy functional, step functions appear as spatial gradients of ReLU functions. Differentiating these step functions with respect to the network parameters, as required when computing gradients during training, yields Dirac's delta distributions concentrated on the ReLU breaking hyperplanes. In the 1-D case, let $b\in (0, 1)$ and consider the following parameter-dependent integral
\[
    j(b) = \int_0^1 \bigl(\mathrm{H}(x - b)\bigr)^2 \,\mathrm{d}x 
         = \int_b^1 \mathrm{d}x = 1-b,
\]
where $\mathrm{H}(\cdot)$ denotes the Heaviside step function arising from the derivative of the $\text{ReLU}$. Differentiating with respect to $b$ yields $j'(b) = -1$. Standard Monte Carlo, however, approximates $j'(b)$ by evaluating 
$\partial_b\bigl(\mathrm{H}(x_i - b)\bigr)^2 = -2\mathrm{H}(x_i-b)\,\delta(x_i - b)$ at 
randomly sampled points. Since $\delta(x_i - b) = 0$ almost everywhere, this gives the estimate $0$ for $j'(b)$, which is incorrect. Consequently, Monte Carlo sampling yields an incorrect gradient regardless of the number of sample points. 
Mathematically, the non-smoothness of ReLU activations prevents the derivative of the energy functional from satisfying the Leibniz integral rule. Therefore, it prevents the interchange of the derivative with respect to the trainable parameters and the limit over the quadrature points: \textbf{The parameter gradient and the integration limit fail to commute}. 
As a result, the training process computes an incorrect gradient and fails to converge to the correct solution. This failure was recently identified in ~\cite{rowan2025failure} for the Deep Ritz Method with single-hidden-layer ReLU networks.

To address this failure in gradient computation, several strategies have been proposed in the literature. One approach employs smooth surrogate gradients during backpropagation \cite{neftci2019surrogate,yin2021accurate, zenke2021remarkable}. These methods use smooth mathematical functions to approximate the Dirac's delta distributions at ReLU breaking hyperplanes. However, this requires setting a smoothing parameter to recover accurate gradients without causing numerical instability. Achieving this balance is difficult in practice, often leading to noisy gradients and poor optimization. Alternatively, others attempt to circumvent the problem by replacing the standard ReLU with smoother variants, such as $\text{ReLU}^2$ or $\text{ReLU}^3$ \cite{so2025higher, yao2026feature, jiao2021rate}. While these higher-order functions eliminate the Dirac's delta distributions arising at the breaking hyperplanes and allow standard stochastic methods to capture the gradient of the energy functional, they introduce some new numerical challenges. 
For instance, variational formulations with ReLU$^k$ networks lead to ill-conditioned linear systems~\cite{mao2026solving}.

Outside the PDE context, a similar problem arises in differentiable rendering, where one integrates the parameter gradient of an image. Since an image is a discontinuous function, its gradient contains Dirac's delta distributions concentrated on the discontinuities~\cite{flanders1973differentiation}. Existing methods sample these contributions on discontinuity edges found by ray tracing~\cite{li2018differentiable}, remove the moving discontinuities by a change of variables~\cite{loubet2019reparameterizing}, or differentiate them systematically~\cite{bangaru2021systematically}. To the best of our knowledge, no existing method retains the standard ReLU activation and integrates both the Deep Ritz energy and its gradients analytically in dimensions higher than one. Consequently, developing a reliable methodology for this setting remains an open research problem.

This paper introduces an integration approach for single-hidden-layer ReLU networks that computes correct gradients when using Deep Ritz for PDEs or PINNs for first-order problems. %
In particular, our approach explicitly identifies the ReLU breaking hyperplanes where Dirac's delta distributions are concentrated. By finding these hyperplanes exactly, our method correctly integrates the gradients, ensuring the optimizer receives correct information to enable convergent training.

Our approach relies on the structure of single-hidden-layer ReLU networks. These networks generate breaking hyperplanes that partition the computational domain into convex polytopes.
For each neuron, we define an \textit{active sub-polytope} as a region where the input to the ReLU is strictly positive, reducing the ReLU to a linear function.
To identify these active sub-polytopes, we utilize the Sutherland–Hodgman algorithm \cite{sutherland1974reentrant}. We then perform integration over each determined active subpolytope. To this end, we introduce exact quadrature via triangulation, which decomposes each active subpolytope into triangles. Because the ReLU function behaves as a linear polynomial within these regions, we apply Gaussian quadrature to evaluate the energy functional and its gradients with machine-precision accuracy. Finally, by assembling the local contributions of each active sub-polytope, we reduce the continuous energy minimization problem to the assembly of a global matrix-vector system, which can be efficiently computed using standard numerical linear algebra methods.

The exact quadrature above relies on the integrand being a polynomial on each active sub-polytope. This holds for the stiffness terms, but not in general for the source term $f$ or for variable PDE coefficients. For such data, exact quadrature is unavailable, and high-order Gaussian rules become expensive and still approximate. To handle these cases, we also propose Conforming Monte Carlo integration, which restricts Monte Carlo sampling exclusively to the identified active sub-polytopes. By confining the sampling to active regions and accounting for their exact boundaries, this approach correctly integrates the Dirac's delta distributions concentrated on the breaking hyperplanes, which ensures gradients are correct up to usual quadrature precision. Hence, Conforming Monte Carlo keeps the flexibility of standard Monte Carlo, while fixing its incorrect gradients.

Finally, to manage the computational complexity of our proposed method, we introduce two training strategies: Neuron-Wise Stochastic and Double Neuron-Wise Stochastic. The stochastic versions essentially sample based on the neurons, and thereby only consider the contributions of selected neurons to the loss function at each iteration. Each strategy offers a different trade-off between accuracy and complexity, allowing users to select the most suitable option depending on the requirements of their specific problem.

While the proposed approach is demonstrated here for elliptic PDEs within the Deep Ritz method using a
single-hidden-layer network, the framework can extend in principle to broader variational PDE formulations involving ReLU networks. The current study focuses on two-dimensional polytopal domains. While the standard Sutherland-Hodgman algorithm admits a formal generalization to arbitrary spatial dimensions, practical deployment in three or higher dimensions faces increased combinatorial complexity. Moreover, for non-polytopal domains, one must resort to curvilinear clipping algorithms or polytopal approximations, thereby incurring a geometric error.

We note that the method proposed in this work is useful beyond the solution of PDEs. Indeed, the same problem arises when networks with discontinuous activations, such as Spiking Neural Networks, are trained with a standard $L^2$-regression loss \cite{neftci2019surrogate,yin2021accurate,zenke2021remarkable}. These works usually rely on surrogate gradients. In contrast, our proposed approach identifies the Dirac's delta contributions, locates the discontinuities, and integrates them exactly.

The remainder of this paper is organized as follows. Section \ref{sec:DRM_ReLU} introduces the Deep Ritz Method and the employed NN architecture. In Section \ref{sec:Limitations}, we analyze the limitations of existing stochastic integration methods. Section \ref{sec:Methodology} details our proposed exact integration framework. This includes the application of the Sutherland-Hodgman algorithm, the matrix-vector energy formulation, and our specialized exact integration and Conforming Monte Carlo techniques. We then present three training algorithms in Section \ref{sec:Training Strategies}, followed by a comprehensive numerical assessment in Section \ref{sec:Numerical experiments}. Finally, we provide concluding remarks and discuss potential extensions in Section \ref{sec:Conclusion}.

\section{Deep Ritz method and ReLU network}\label{sec:DRM_ReLU}
In this section, we describe the model elliptic problem and the Deep Ritz method. Then, we introduce a single-hidden-layer network used as the trial space for the associated minimization problem in the Deep Ritz framework.

\subsection{Model problem} %
Let $\Omega \subset \mathbb{R}^d$, with $d \in \{1, 2\}$, denote a bounded convex polytope, and let $L^2(\Omega)$ be the space of square-integrable functions on $\Omega$. For any integer $k \ge 1$, we denote by $H^k(\Omega)$ the Sobolev space of functions in $L^2(\Omega)$ whose weak derivatives up to order $k$ also belong to $L^2(\Omega)$. Let $H_0^1(\Omega) = \{v \in H^1(\Omega): v|_{\partial\Omega} = 0\}$ be the standard Sobolev space of order one with zero trace on the boundary. We consider the following elliptic boundary value problem: find $u \in H^1_0(\Omega)$ such that 
\begin{equation}\label{eq_PDE}
\left\{
\begin{array}{r@{}ll}
 - \Delta u  + a \,u &{}= f, 
  & \text{in }\Omega\,,
 \\
 u &{}= 0, & \text{on }\partial\Omega\,,
\end{array}\right.
\end{equation}
where the equation is understood in the weak sense, $a \geq 0$ is a constant reaction coefficient, and $f \in L^2(\Omega)$ represents the source term. Equivalently, this problem can be characterized through a variational energy framework~\cite[Section 2.2.5]{evans2010partial}. The solution $u$ is the unique minimizer of the energy functional $\mathcal{J}: H^1_0(\Omega) \rightarrow \mathbb{R}$, defined by
\begin{equation} \label{energy_function}
\mathcal{J}(v) \coloneqq \int_\Omega \left( \frac{1}{2} |\nabla v|^2 + \frac{1}{2} a v^2 - f v \right) dx.
\end{equation}
The corresponding variational problem then reads
\begin{equation}
\left\{
\begin{array}{l}
   \text{find } u \in H^1_0(\Omega) \text{ such that:}   \\[2pt]
   \mathcal{J}(u) \le \mathcal{J}(v),  \quad \forall\, v \in H^1_0(\Omega).
\end{array}
\right.
\end{equation}

\subsection{The Deep Ritz method}
The Deep Ritz Method \cite{yu2018deep} approximates the solution of the elliptic PDE~\eqref{eq_PDE} by minimizing the corresponding energy functional~\eqref{energy_function} over a trial space defined by an NN. To define it, we first introduce an architecture given by a single-hidden-layer network $\mathcal{N}_{\bstheta}:= \mathcal{N}_{n,\bstheta}: \mathbb{R}^d \rightarrow \mathbb{R}$ with $n$ neurons, i.e.,
\begin{equation}\label{NN-definition}
\mathcal{N}_{\bstheta}(\bs{x})= \sum_{j=1}^n c_j \ReLU( \bs{w}_j \cdot \bs{x} + b_j) 
= \bsc \cdot \ReLU (W \bs{x} + \bsb),
\end{equation}
where $\bstheta$ denotes the collection of all parameters $(\bs{c}, W, \bs{b})$ in the network. Here $\bs{c} \in \mathbb{R}^n$ represents the outer weights, $W \in \mathbb{R}^{n \times d}$ is the weight matrix (with rows $\bs{w}_j$), and $\bs{b} \in \mathbb{R}^n$ is the bias vector. 

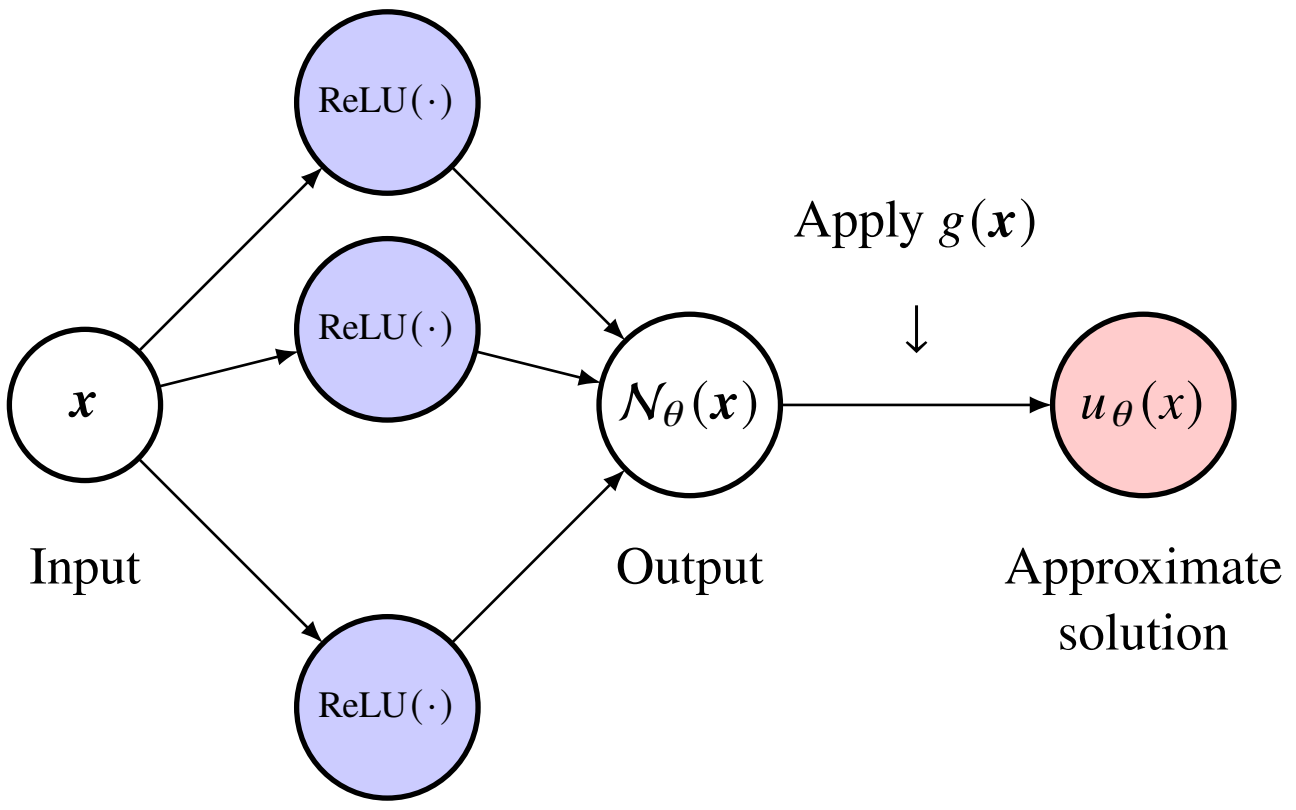
\begin{figure}[htbp]
	\begin{center}
	\input{Figure_Base.tex}
	\end{center}
	\caption{Boundary-conforming single-hidden-layer NN architecture.}
	\label{fig_NN_architecture_u}
\end{figure}

For a correct approximation of the solution to problem~\eqref{eq_PDE} via NNs, the approximating function must also satisfy the boundary conditions. This can usually be done in two ways, via a strong or a weak imposition \cite{berg2018unified}. 
In our case, we strongly enforce the essential homogeneous boundary conditions by defining a trial solution that follows a conforming structure \cite{berrone2023enforcing}, i.e.,
\begin{equation}\label{eq_DBC_NN}
u(\bs{x}) \approx u_{\bstheta}(\bs{x}) = \mathcal{N}_{\bstheta}(\bs{x})g(\bsx).
\end{equation}
Herein, $g$ is a bubble function that vanishes on $\partial \Omega$ and remains positive within $\Omega$. In our case, $g(\bsx)$ is constructed as a product of algebraic boundary equations.
This construction ensures that $u_{\bstheta}(\bs{x})$ strongly satisfies zero boundary conditions, allowing the network to focus on learning the interior behavior of the solution. Given that any nonhomogeneous PDE problem can be transformed into a homogeneous one, this methodology can be readily adapted to handle inhomogeneous boundary conditions by applying an appropriate lifting %
function. The schematic representation of this boundary-conforming NN architecture is illustrated in Figure \ref{fig_NN_architecture_u}.

The set of boundary-conforming NN functions $u_{\bstheta}$ defined in \eqref{eq_DBC_NN} constitutes our approximation space $\mathbb{M}_n$, which serves as a trial space within the Sobolev space $H^1_0(\Omega)$. Specifically, we define this space as  
\begin{alignat}{2}
    \mathbb{M}_n \coloneqq \left\{
    u_{\bstheta}\in H^1_0(\Omega) \;\big|\;
    \, \bstheta\in \bsTheta
    \right\},
\end{alignat}
where the parameter space $\bsTheta:= \mathbb{R}^{n} \times \mathbb{R}^{n \times d} \times \mathbb{R}^{n} \simeq \mathbb{R}^{n \times (d+2)}$ is defined by the dimensions of the weights and biases. 
The Deep Ritz Method seeks the optimal parameters $\bstheta^*$ by minimizing the variational energy functional $\mathcal{J}(\cdot)$ over the set $\mathbb{M}_n$, i.e.,
\begin{equation}\label{eq:DR-minimization}
\left\{
\begin{array}{l}
   \text{find }  \bstheta^* \in \bsTheta \text{ such that:}   \\[2pt]
   \mathcal{J}(u_{\bstheta^*}) \le \mathcal{J}(u_{\bstheta}),  \quad \forall\,\bstheta \in \bsTheta.
\end{array}
\right.
\end{equation}
Unlike traditional Ritz methods that use linear subspaces (e.g., finite element spaces), the set $\mathbb{M}_n$ is nonlinear, thereby transforming the problem into a nonconvex optimization problem typically solved via gradient-based methods; see, e.g.,~\cite{nocedal_numerical_2006}.
\section{Monte Carlo Integration and Surrogate Gradient: Gradient Bias from ReLU Activations}\label{sec:Limitations}

In practical implementations of the Deep Ritz Method, training the NN to find the optimal parameters $\bstheta^*$ involves minimizing the variational energy functional~\eqref{energy_function}. 
To evaluate this functional, one typically relies on numerical sampling and updates the network parameters using gradient-based optimization algorithms, such as Stochastic Gradient Descent (SGD) or Adam~\cite{AdamOpt}. Although these frameworks are highly effective for smooth activation functions, they face a numerical integration problem when applied to networks with ReLU activations.

This section examines the core difficulties associated with the standard Monte Carlo integration method and its heuristic approaches, such as using a surrogate gradient. We highlight why these approaches produce incorrect gradient computations, thereby establishing the primary motivation for our proposed approach.

\subsection{Standard Monte Carlo}
Let $\{\bs{x}_i\}_{i=1}^N$ be a set of $N$ independent and identically distributed (i.i.d.) points sampled from a uniform distribution over $\Omega$. The standard Monte Carlo approximation for the energy functional \eqref{energy_function} is defined as follows
\begin{equation} \label{eq_Monte Carlo_functional}
\mathcal{J}_{_{MC}}({\bstheta}) \coloneqq \frac{|\Omega|}{N} \sum_{i=1}^N \left( \frac{1}{2} |\nabla u_{\bstheta}(\bs{x}_i)|^2 + \frac{1}{2} a u_{\bstheta}^2(\bs{x}_i) - f(\bs{x}_i) u_{\bstheta}(\bs{x}_i) \right),
\end{equation}
where $|\Omega|$ denotes the volume of the domain $\Omega$ and $\{\bs{x}_i\}$ is resampled at each training iteration. This estimator is unbiased when the integrand is sufficiently regular, converging to the exact integral as $N \to \infty$. 
However, when $u_{\bstheta}$ uses ReLU activations, the Monte Carlo gradient estimator is biased. Because the non-smoothness of ReLU activations violates the Leibniz integral rule, the parameter derivative and the integration limit fail to commute. Standard Monte Carlo implicitly interchanges these non-commutative limits and evaluates the parameter derivative pointwise, which vanishes almost everywhere. This misses the Dirac’s delta distributions produced on the breaking hyperplanes and yields incorrect gradients. 

To illustrate this bias, we revisit the one-dimensional example of the introduction, $j(b) = \int_0^1 \bigl(\mathrm{H}(x-b)\bigr)^2 \, dx = 1-b$, whose exact derivative is $j'(b) = -1$. The integrand equals one on $(b,1)$ and zero on $(0,b)$, and it does not depend on $b$ inside either interval. Hence, the whole derivative comes from the moving point $x=b$: increasing $b$ shrinks the interval $(b,1)$. Standard Monte Carlo approximates $j(b)$ by
\begin{equation}\label{eq:jMC}
    \mathcal{J}_{MC}(b) = \frac{1}{N}\sum_{i=1}^{N} \bigl(\mathrm{H}(x_i-b)\bigr)^2 ,
\end{equation}
which converges to $j(b)$ as $N \to \infty$. However, its derivative is evaluated pointwise, i.e., $\mathcal{J}_{MC}'(b) = -\frac{2}{N}\sum_{i=1}^{N} \mathrm{H}(x_i-b)\,\delta(x_i-b) = 0$, since, with probability one, no sample lies exactly at $x=b$. Therefore,
\begin{equation}\label{eq:limit_commute}
    \lim_{N\to\infty} \mathcal{J}_{MC}'(b) = 0 \;\neq\; -1
    = \frac{d}{db}\,\lim_{N\to\infty} \mathcal{J}_{MC}(b).
\end{equation}
The missing value is exactly the Dirac's delta contribution, $-2\int_0^1 \mathrm{H}(x-b)\,\delta(x-b)\,dx = -2\,\mathrm{H}(0) = -1$, where we use the convention $\mathrm{H}(0) = 1/2$. In the Deep Ritz energy \eqref{energy_function}, the same contributions are concentrated on the breaking hyperplanes, which have zero Lebesgue measure. Hence, standard Monte Carlo misses them for any $N$.

\subsection{Using Surrogate Gradients}

To address the limitations of standard Monte Carlo integration in capturing the Dirac's delta contributions arising at the breaking hyperplanes, an alternative is to employ surrogate gradients \cite{neftci2019surrogate,yin2021accurate, zenke2021remarkable}. Under this approach, the neural network forward pass proceeds as usual, using the standard $\text{ReLU}$ activation function. However, during the backward pass, the gradient with respect to the network parameters, which gives rise to a Dirac's delta distribution, is replaced with a smooth analytical function that approximates this singular contribution.

In this context, a Gaussian surrogate is employed to replace the Dirac’s delta distribution
\cite{yin2021accurate}
\begin{equation}\label{Gaussian surrogate gradients}
    \delta_\beta(x) = \frac{1}{\sqrt{2\pi\beta}}\exp\left(-\frac{x^2}{2\beta}\right),
\end{equation}
where $\beta > 0$ is a smoothing parameter that controls the sharpness of the approximation; as $\beta\rightarrow 0$, $\delta_\beta$ converges to the Dirac's delta distribution in the distributional sense. This provides a computationally tractable approach that enables gradient-based optimizers to compute approximate gradients at the non-differentiable activation boundaries, allowing network parameters to receive gradient updates during training.

\noindent However, this approach introduces some limitations. The selection of the smoothing parameter $\beta$ presents a strict trade-off between the accuracy of the approximation and the stability of the Monte Carlo sampling. Minimizing the approximation error requires $\beta \to 0$ so that $\delta_\beta$ converges to the true Dirac's delta distribution. However, as $\beta \to 0$, the peak height of the Gaussian surrogate gradient \eqref{Gaussian surrogate gradients} grows as $\mathcal{O}(\beta^{-1/2})$ while the effective support shrinks rapidly. As a result, randomly placed Monte Carlo points almost always miss this narrow peak. Consequently, most sampled points evaluate to near zero, while a rare point landing near the center yields an extremely large value of order $\mathcal{O}(\beta^{-1/2})$. This extreme gap between near-zero and very large values across different Monte Carlo batches causes the sample variance of the gradient estimator to explode \cite{mohamed2020monte}. In practice, this creates two destructive failure modes during training: when a batch misses the peak entirely, the computed gradient vanishes and parameter updates freeze; when a batch happens to hit the peak, it injects an excessively large gradient step that causes training to oscillate wildly. Consequently, maintaining stability under small $\beta$ demands costly computational workarounds, such as dynamically increasing the sample size $N$ while aggressively damping the learning rate throughout training. Balancing these competing demands, namely narrowing the gradient approximation while simultaneously scaling the Monte Carlo sample size and learning rate, is challenging. This highlights the critical need for a framework capable of computing the exact singular Dirac's delta contributions arising at the breaking hyperplanes explicitly.

\section{Exact Integration Methodology}\label{sec:Methodology}
This section introduces two integration methodologies that overcome the aforementioned limitations, enabling an accurate computation of the gradient in problem~\eqref{eq:DR-minimization}. 

We first differentiate the trial solution, yielding:
\begin{equation}\label{Gradient_solution_1}
\nabla_{\bsx}\, u_{\bstheta}(\bsx) = \nabla_{\bsx}\, \mathcal{N}_{\bstheta}(\bsx)g(\bsx) + \mathcal{N}_{\bstheta}(\bsx) \nabla_{\bsx}\, g(\bsx).
\end{equation}
Furthermore, the gradient of the network with respect to the spatial coordinates $\bsx$ is given by 
\begin{equation}\label{Gradient_N_1}
    \nabla_{\bsx}\, \mathcal{N}_{\bstheta}(\bsx) = \sum_{i=1}^n c_i\bsw_i \, \text{H}(\bsw_i \cdot \bsx + b_i).
\end{equation}
Substituting \eqref{Gradient_N_1} into \eqref{Gradient_solution_1}, we obtain the following explicit form for the gradient of the approximate solution
\begin{equation}\label{Gradient_solution1}
\nabla_{\bsx}\, 
u_{\bstheta}(\bsx) = \left(\sum_{i=1}^n c_i\bsw_i \, \text{H}(\bsw_i \cdot  \bsx + b_i)\right) g(\bsx) + \left( \sum_{i=1}^n c_i \, \text{ReLU}(\bsw_i \cdot \bsx + b_i) \right) \nabla_{\bsx}\, g(\bsx).
\end{equation}
Consequently, the energy functional 
$\mathcal{J}(u_{\bstheta})$ can be decomposed into three terms 
\begin{equation}\label{J_I_II_III}
\mathcal{J}(u_{\bstheta}) = \frac{1}{2} \mathrm{I}(\bstheta) + \frac{1}{2}\mathrm{II}(\bstheta) - \mathrm{III}(\bstheta),
\end{equation}
where
\begin{equation}\label{Is1}
\begin{array}{rl}
\mathrm{I}(\bstheta) =  & \hspace{-0.2cm} \displaystyle \int_\Omega |\nabla u_{\bstheta}(\bsx)|^2 d\bsx = \displaystyle \int_\Omega \left| \left(\sum_{i=1}^n c_i\bsw_i  \text{H}(\bsw_i \cdot \bsx + b_i)\right) g(\bsx) + \left( \sum_{i=1}^n c_i \text{ReLU}(\bsw_i \cdot \bsx + b_i) \right) \nabla g(\bsx) \right|^2 \, d\bsx, \\[4mm]
\mathrm{II}(\bstheta) = & \hspace{-0.2cm} \displaystyle \int_\Omega a (u_{\bstheta}(\bsx))^2 \, d\bsx = a \int_\Omega  \left(\sum_{i=1}^n c_i \text{ReLU}(\bsw_i \cdot \bsx + b_i)g(\bsx) \right)^2 \, d\bsx, \\[4mm]
\mathrm{III}(\bstheta) = & \hspace{-0.2cm} \displaystyle \int_\Omega f(\bsx) u_{\bstheta}(\bsx) \, d\bsx = \int_\Omega f(\bsx) \left( \sum_{i=1}^n c_i \text{ReLU}(\bsw_i \cdot \bsx + b_i)g(\bsx) \right) \, d\bsx. \
\end{array}\end{equation}

In the following subsections, we describe the Sutherland–Hodgman algorithm~\cite{sutherland1974reentrant} to partition $\Omega$ into active sub-polytopes, explicitly capturing the activation interfaces; next, we derive a compact matrix-vector representation of $\mathcal{J}(u_{\bstheta})$; we then introduce two integration methods, an exact integration and a Conforming Monte Carlo, to evaluate the integrals over active subdomains; and finally, we analyze the computational cost of the geometric preprocessing step.

\subsection{Clipping polygons: Sutherland-Hodgman algorithm}
To evaluate \eqref{J_I_II_III} and \eqref{Is1}, we must precisely determine the intersection of the convex polytope $\Omega$ with the support of the activation function and its gradient. Specifically, the functions $\text{ReLU}(\boldsymbol{w}_i \cdot \boldsymbol{x} + b_i)$ and $\text{H}(\boldsymbol{w}_i \cdot \boldsymbol{x} + b_i)$ are supported on the half-space $\mbf{X}_i^+ :=\{ \boldsymbol{x} \in \mathbb{R}^d \mid \boldsymbol{w}_i \cdot \boldsymbol{x} + b_i \ge 0 \}$, for $i = 1, \ldots, n$. Consequently, for each neuron $i$, we define the active subdomain $\Omega_i$ as the intersection of the physical domain $\Omega$ and the half-space $\mbf{X}_i^+$, i.e.,
\begin{equation}
\Omega_i = \Omega \cap \mbf{X}_i^+ = \{ \boldsymbol{x} \in \Omega \mid 
\boldsymbol{w}_i \cdot \boldsymbol{x} + b_i \ge 0 \}, \qquad i = 1, \ldots, n.
\end{equation}
Since $\Omega$ is a convex polytope and $\mbf{X}_i^+$ is a half-space, $\Omega_i$ is a convex polytope. To construct $\Omega_i$ computationally, we employ the Sutherland--Hodgman algorithm \cite{sutherland1974reentrant} to perform polygon clipping, that is, computing the exact intersection between the domain polygon $\Omega$ and the half-space $\mathbf{X}_i^+$. 
Algorithm~\ref{Sutherland--Hodgman Algorithm} below describes this 
procedure: it processes the edges of $\Omega$ sequentially against the 
bounding hyperplane of $\mathbf{X}_i^+$, retaining vertices inside the 
half-space and inserting new vertices at boundary intersection points.

\begin{algorithm}[H]
{
\begingroup
\SetAlgoLined
Select the ordered vertices of $\Omega$, i.e., $\mathcal{P} = \{\bsp_1, \ldots, \bsp_m\}$ 
and the hyperplane coefficients $(\boldsymbol{w}_i, b_i)$ of $\mathbf{X}_i^+$. 
Initialize an empty output list $\mathcal{Q}$.\\

\For{$k = 1, \ldots, m$ (indices taken modulo $m$)}{
Select two consecutive vertices $\bsp_k$ and $\bsp_{k+1}$;\\
Compute the signed distances $\bsw_i \cdot \bsp_k + b_i$ and 
$\bsw_i \cdot \bsp_{k+1} + b_i$ to the hyperplane;\\
\If{ the endpoints lie on opposite sides}{
\hspace{0.5cm} Compute the intersection point 
by linear interpolation and append it to $\mathcal{Q}$;}
\If{$\bsp_{k+1}$ lies inside $\mathbf{X}_i^+$}{Append $\bsp_{k+1}$ to $\mathcal{Q}$;}
}
Return $\mathcal{Q}$ as the ordered vertices of $\Omega_i = \Omega \cap \mathbf{X}_i^+$.
\caption{Sutherland--Hodgman polygon clipping}\label{Sutherland--Hodgman Algorithm}
\endgroup
}
\end{algorithm}

The following proposition characterizes the vertices returned by 
Algorithm~\ref{Sutherland--Hodgman Algorithm} explicitly in terms of the 
network parameters $\bstheta$.

\begin{proposition}\label{proposition_clipping}
Let $\Omega$ be a convex polygon defined by its ordered vertices
$\mathcal{P} = \{\bsp_1, \ldots, \bsp_m\}$,
and let $\Gamma_i := \{\bsx : \bsw_i \cdot \bsx + b_i = 0\}$ denote the
activation hyperplane of the $i$-th neuron.
Assume $\operatorname{int}(\Omega) \cap \Gamma_i \neq \emptyset$.
Then $\bsp$ is a vertex of $\Omega_i$ if
\begin{equation}\label{eq_proposition}
\left\{
\begin{alignedat}{2}
  &\bsw_i \cdot \bsp + b_i > 0, \\
  &\bsp \in \mathcal{P},
\end{alignedat}
\right.
\qquad \text{or} \qquad
\left\{
\begin{alignedat}{2}
  &\bsw_i \cdot \bsp + b_i = 0, \\
  &\bsp \in \operatorname{Conv}(\bsp_k,\bsp_{k+1}), \quad \text{for some } k,
\end{alignedat}
\right.
\end{equation}
where $\operatorname{Conv}(\bsp_k,\bsp_{k+1})$ denotes the closed line segment
between $\bsp_k$ and $\bsp_{k+1}$.
\end{proposition}
\begin{proof}
See Appendix~\ref{Proof of Proposition}.
\end{proof}

\begin{corollary}\label{corollary_differentiability}
Under the assumptions of Proposition~\ref{proposition_clipping}, the vertices 
of $\Omega_i$ are differentiable almost everywhere with respect to the network parameters 
$(\bsw_i, b_i)$. Vertices satisfying the second condition in \eqref{eq_proposition} are of the form $\bsp = \bsp_k + t\,(\bsp_{k+1} - \bsp_k)$, where $t \in [0,1]$ is the unique parameter such that $\bsp \in \Gamma_i$, given by
\begin{equation}
t = \frac{-(\bsw_i \cdot \bsp_k + b_i)}{\bsw_i \cdot (\bsp_{k+1} - \bsp_k)}.
\end{equation}
Since $t$ is a smooth function of $(\bsw_i, b_i)$, each intersection vertex depends smoothly on the network parameters. Consequently, the active subdomain $\Omega_i$ and all integrals over it are differentiable almost everywhere with respect to $\bstheta$, enabling gradient computation via automatic differentiation.
\end{corollary}

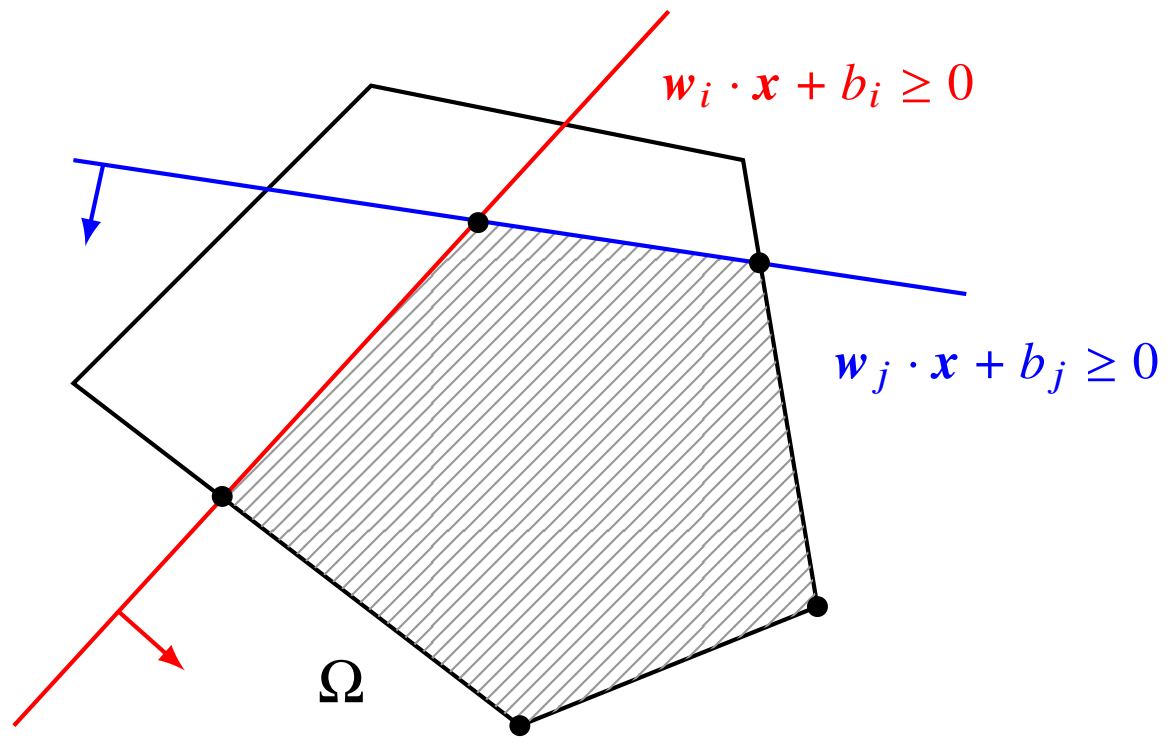
\begin{figure}[htbp]
    \centering
    \resizebox{0.5\textwidth}{!}{ %
        \input{Clipping_polygons.tex}
    }
    \caption{Geometric clipping of the domain $\Omega$ via Proposition \ref{proposition_clipping}.}
    \label{fig_Geometric clipping}
\end{figure}

\noindent Evaluation of the energy functional $\mathcal{J}(u_{\bstheta})$ 
requires integrating the products of the ReLU and Heaviside functions. 
Therefore, for a pair of neurons $i$ and $j$, the product is non-zero only 
on the intersection $\Omega_{ij} = \Omega \cap \mbf{X}_i^+ \cap \mbf{X}_j^+$. 
As illustrated in Figure~\ref{fig_Geometric clipping}, 
Proposition~\ref{proposition_clipping} can be applied recursively. Indeed, 
to find the region where multiple neurons are simultaneously ``active'', we 
clip $\Omega$ against the first hyperplane to obtain $\Omega_i$, then clip 
$\Omega_i$ against the second hyperplane to obtain $\Omega_{ij}$. This 
ensures that the complex global integral is reduced to a sum of exact 
integrals over well-defined convex polygons.

\subsection{Matrix-vector formulation of the energy functional}
By applying Proposition \ref{proposition_clipping} to partition the domain into active sub-polytopes, we can express the integrands in Equation \eqref{Is1} as piecewise polynomials. Therefore, the components of the energy functional are as follows
\begin{equation}\label{Is_homogeneous1}
\begin{array}{rl}
\mathrm{I}(\bstheta) =  & \displaystyle \sum_{i = 1}^n\sum_{j =1 }^n\int_{\Omega_{ij}} \Big(  c_i\bsw_i g(\bsx) + c_i (\bsw_i \cdot \bsx + b_i) \nabla g \Big) \cdot \Big( c_j\bsw_j g(\bsx) + c_j (\bsw_j \cdot \bsx + b_j) \nabla g  \Big) \, d\bsx, \\[4mm]
\mathrm{II}(\bstheta) = & \displaystyle \sum_{i = 1}^n\sum_{j =1 }^n a \int_{\Omega_{ij}}  \Big( c_i(\bsw_i \cdot \bsx + b_i)g(\bsx) \Big)\Big( c_j(\bsw_j \cdot \bsx + b_j)g(\bsx) \Big) \, d\bsx, \\[4mm]
\mathrm{III}(\bstheta) = & \displaystyle \sum_{i = 1}^n\int_{\Omega_{i}}  \Big( c_i(\bsw_i \cdot \bsx + b_i)g(\bsx)f(\bsx)\Big) \, d\bsx. \
\end{array}\end{equation}

The following propositions establish the transformation of the integrals \eqref{Is_homogeneous1} into a compact matrix-vector representation. 
\begin{proposition}\label{proposition_I}
Let $\bsz_i = c_i \bsw_i \in \mathbb{R}^{d\times 1}$ and $y_i = c_i b_i\in \mathbb{R}$. The integral $\mathrm{I}(\bstheta)$ admits the following  matrix representation 
\begin{equation}\label{Matrix-form-I}
\mathrm{I}(\bstheta) = \bsz^\top \mathcal{A}^{(1)} \bsz + 2 \, \bsz^\top \mathcal{A}^{(2)} \bsy + \bsy^\top \mathcal{A}^{(3)} \bsy,
\end{equation}
where $\bsz = [\bsz_i] \in \mathbb{R}^{nd\times1}$ and $\bsy = [y_i]\in \mathbb{R}^{n\times1}$. Furthermore, the global operator matrices $\mathcal{A}^{(1)} =[\mathcal{A}^{(1)}_{ij}] \in \mathbb{R}^{nd \times nd}$, $\mathcal{A}^{(2)}= [\mathcal{A}^{(2)}_{ij}] \in \mathbb{R}^{nd \times n}$, and $\mathcal{A}^{(3)} = [\mathcal{A}^{(3)}_{ij}]\in \mathbb{R}^{n \times n}$ are assembled from pairwise neuron interactions over the clipped domains $\Omega_{ij}$ and are defined as follows
\begin{equation}\label{defining_A123}
\mathcal{A}^{(1)}_{ij} = \int_{\Omega_{ij}} \mathcal{M}^\top \mathcal{M} \, d\bsx, \quad
\mathcal{A}^{(2)}_{ij} = \int_{\Omega_{ij}} \mathcal{M}^\top \nabla g \, d\bsx, \quad
\mathcal{A}^{(3)}_{ij} = \int_{\Omega_{ij}} |\nabla g |^2 \, d\bsx,
\end{equation}
where $\mathcal{M}(\bsx) = g(\bsx) \mathbf{I} + \nabla g(\bsx) \bsx^\top \in \mathbb{R}^{d \times d}$ with the $d \times d$ identity matrix $\mathbf{I}$.
\end{proposition}
\begin{proof}
Substituting the auxiliary variables $\bsz_i$ and $y_i$ into the definition of $\mathrm{I}(\bstheta)$, we obtain
\begin{equation}
\mathrm{I}(\bstheta) = \sum_{i = 1}^n\sum_{j =1 }^n \int_{\Omega_{ij}} \left( \mathcal{M}(\bsx) \bsz_i + y_i \nabla g(\bsx) \right)^\top \left( \mathcal{M}(\bsx) \bsz_j + y_j \nabla g(\bsx) \right) \, d\bsx.
\end{equation}
Expanding the integrand yields
\begin{equation}\label{I_1}
\begin{array}{rl}
\displaystyle \mathrm{I}(\bstheta) = \sum_{i= 1}^n\sum_{j = 1}^n \Big( &\displaystyle \hspace{-0.2cm}\bsz_i^\top \left[ \int_{\Omega_{ij}} \mathcal{M}^\top \mathcal{M} \,d\bsx \right] \bsz_j + \bsz_i^\top \left[ \int_{\Omega_{ij}} \mathcal{M}^\top \nabla g \, d\bsx \right] y_j +\\
&\displaystyle  y_i \left[ \int_{\Omega_{ij}} \nabla g^\top \mathcal{M} \, d\bsx \right] \bsz_j + y_i \left[ \int_{\Omega_{ij}} |\nabla g |^2 \, d\bsx \right] y_j \Big)
\end{array}
\end{equation}
According to \eqref{I_1}, the precomputable local kernels \eqref{defining_A123} lead directly to the block matrix representation given in \eqref{Matrix-form-I}.
\end{proof}
\begin{proposition}\label{proposition_II}
The reaction term $\mathrm{II}(\bstheta)$ is formulated in the following matrix form
\begin{equation}\label{Matrix-form-II}
\mathrm{II}(\bstheta) = a \left( \bsz^\top \mathcal{B}^{(1)} \bsz + 2 \,\bsz^\top \mathcal{B}^{(2)} \bsy + \bsy^\top \mathcal{B}^{(3)} \bsy \right).
\end{equation}
The global kernel matrices $\mathcal{B}^{(1)} = [\mathcal{B}^{(1)}_{ij}] \in \mathbb{R}^{nd \times nd}$, $\mathcal{B}^{(2)} = [\mathcal{B}^{(2)}_{ij}]\in \mathbb{R}^{nd \times n}$, and $\mathcal{B}^{(3)} = [\mathcal{B}^{(3)}_{ij}]\in \mathbb{R}^{n \times n}$ are assembled from pairwise neuron interactions over the clipped domains $\Omega_{ij}$ as follows
\begin{equation}\label{defining_B123}
\mathcal{B}^{(1)}_{ij} = \int_{\Omega_{ij}} \mathcal{D} \mathcal{D}^\top \, d\bsx, \quad
\mathcal{B}^{(2)}_{ij} = \int_{\Omega_{ij}} \mathcal{D} g \, d\bsx, \quad
\mathcal{B}^{(3)}_{ij} = \int_{\Omega_{ij}} g^2 \, d\bsx,
\end{equation}
where $\mathcal{D}(\bsx) = g(\bsx)\bsx \in \mathbb{R}^{d \times 1}$.
\end{proposition}
\begin{proof}
Substituting the auxiliary variables $\bsz_i$ and $y_i$, and $\mathcal{D}(\bsx)$ into the definition of $\mathrm{II}(\bstheta)$, we obtain
\begin{equation}
\mathrm{II}(\bstheta) = a \sum_{i = 1}^n\sum_{j = 1}^n  \int_{\Omega_{ij}} \left( \mathcal{D}(\bsx)^\top \bsz_i + y_i g(\bsx) \right) \left( \mathcal{D}(\bsx)^\top \bsz_j + y_j g(\bsx) \right) \, d\bsx.
\end{equation}
Expanding the product within the integrand results in
\begin{equation}\label{II_1}
\mathrm{II}(\bstheta) = a \sum_{i = 1}^n \sum_{j=1}^n \left( \bsz_i^\top \left[ \int_{\Omega_{ij}} \mathcal{D} \mathcal{D}^\top \, d\bsx \right] \bsz_j + \bsz_i^\top \left[ \int_{\Omega_{ij}} \mathcal{D} g \, d\bsx \right] y_j + y_i \left[ \int_{\Omega_{ij}} g \mathcal{D}^\top \, d\bsx \right] \bsz_j + y_i \left[ \int_{\Omega_{ij}} g^2 \, d\bsx \right] y_j \right).
\end{equation}
According to \eqref{II_1}, the precomputable local kernels \eqref{defining_B123} lead directly to the block matrix representation given in \eqref{Matrix-form-II}.
\end{proof}

\begin{proposition}\label{proposition_III}
The linear source term $\mathrm{III}(\bstheta)$ admits the following matrix representation
\begin{equation}\label{Matrix-form-III}
\mathrm{III}(\bstheta) = \bsz^\top \mathcal{C}^{(1)} + \bsy^\top \mathcal{C}^{(2)}.
\end{equation}
The global load vectors $\mathcal{C}^{(1)} = [\mathcal{C}^{(1)}_i] \in \mathbb{R}^{nd \times 1}$ and $\mathcal{C}^{(2)} = [\mathcal{C}^{(2)}_i] \in \mathbb{R}^{n \times 1}$ are assembled from individual neuron contributions over the active subdomains $\Omega_i$ and are defined as follows
\begin{equation}\label{defining_C12}
\mathcal{C}^{(1)}_i = \int_{\Omega_{i}} \mathcal{H}(\bsx) \, d\bsx \in \mathbb{R}^{d \times 1}, \quad \mathcal{C}^{(2)}_i = \int_{\Omega_{i}} f(\bsx)g(\bsx) \, d\bsx \in \mathbb{R},
\end{equation}
where $\mathcal{H}(\bsx) = f(\bsx)g(\bsx)\bsx \in \mathbb{R}^{d \times 1}$.
\end{proposition}
\begin{proof}
Substituting the auxiliary variables $\bsz_i$ and $y_i$ into the definition of the integral $\mathrm{III}(\bstheta)$, we obtain
\begin{equation}\label{III_1}
\mathrm{III}(\bstheta) = \sum_{i = 1}^n \int_{\Omega_{i}} \left( \bsz_i^\top \mathcal{H}(\bsx) + y_i f(\bsx)g(\bsx) \right) \, d\bsx.
\end{equation}
According to \eqref{III_1}, the precomputable local kernels \eqref{defining_C12} on each active subdomain $\Omega_i$ lead directly to the block matrix representation given in \eqref{Matrix-form-III}.
\end{proof}
\begin{corollary}\label{corollary_matrix_form}
By combining Propositions \ref{proposition_I}--\ref{proposition_III}, the total energy functional $\mathcal{J}(u_{\bstheta})$ is expressed as the following quadratic form
\begin{equation}\label{J_matrix_fom1}
\mathcal{J}(u_{\bstheta}) = \frac{1}{2} \bsr^\top \mathcal{K} \bsr - \bsr^\top \mathcal{F},
\end{equation}
where $\bsr = [\bsz^\top, \bsy^\top]^\top \in \mathbb{R}^{n(d+1) \times 1}$ is the augmented parameter vector. The symmetric global stiffness matrix $\mathcal{K} \in \mathbb{R}^{n(d+1) \times n(d+1)}$ and the global load vector $\mathcal{F} \in \mathbb{R}^{n(d+1) \times 1}$ are defined as
\begin{equation}
\mathcal{K} = \begin{bmatrix}
\mathcal{A}^{(1)} + a\mathcal{B}^{(1)} & \mathcal{A}^{(2)} + a\mathcal{B}^{(2)} \\
\left(\mathcal{A}^{(2)} + a\mathcal{B}^{(2)}\right)^\top & \mathcal{A}^{(3)} + a\mathcal{B}^{(3)}
\end{bmatrix}, \quad
\mathcal{F} = \begin{bmatrix} \mathcal{C}^{(1)} \\ \mathcal{C}^{(2)} \end{bmatrix}.
\end{equation}
\end{corollary}

\subsection{Computing the integral}
To assemble the global system \eqref{J_matrix_fom1}, one must evaluate the domain integrals over the active subdomains. We present two strategies for this purpose. The first is an exact integration scheme based on triangular decomposition, which is applicable when the right-hand side $f$ is a polynomial. The second is a Conforming Monte Carlo method, which applies to general $f$.
 
\subsubsection{Exact Integration via Triangular Decomposition}
Assume that $f$ and $g$ are polynomials, and that the Gauss rule is exact for the resulting polynomial degree. Since the spatial kernels defined over the polygonal subdomains $\Omega_{ij}$ are piecewise polynomial, the corresponding integrals can be evaluated exactly. To this end, we adopt a triangular decomposition strategy, which reduces integration over an arbitrary convex polygon to a sum of integrals over standard triangular elements. Each triangle is then handled via high-order Gaussian quadrature applied through Duffy's transformation \cite{duffy1982quadrature}, yielding machine-precision results independent of the geometric complexity of the polygon. Given an active subdomain $\Omega_{ij}$ defined by the ordered set of $K$ vertices $\{\bs{p}_k\}_{k=1}^K$, the integration is performed according to the following algorithm.

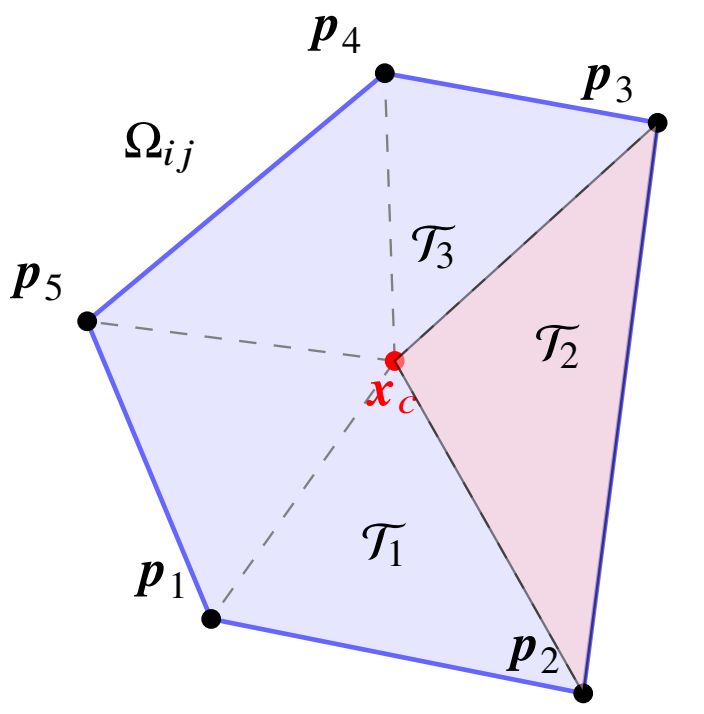
\begin{figure}[htbp]
    \centering
    \resizebox{0.3\textwidth}{!}{ %
        \input{Fig_triangles.tex}
    }
    \caption{Geometric representation of polygon triangulation.}
    \label{fig_triangulation}
\end{figure}

\begin{algorithm}[H]\label{Algorithm_1}
{
\begingroup

\SetAlgoLined
Identify the geometric center $\bsx_c$ of the polygon $\Omega_{ij}$ (or $\Omega_{i}$) by averaging its vertex coordinates $\{\bsp_k\}_{k = 1}^K$.\\
Decompose the polygon into $K$ triangles $\{\mathcal{T}_k\}_{k = 1}^K$ by connecting $\bsx_c$ to each consecutive pair of vertices $(\bsp_k, \bsp_{k+1})$, with indices taken modulo $K$.\\
Map each physical triangle to a standard reference triangle to provide a consistent coordinate system for the quadrature rule.\\
Apply a two-dimensional Legendre-Gauss quadrature rule to each triangle.\\
Sum the results from all triangles to obtain the final integral value for the polygon $\Omega_{ij}$.
\caption{Exact Integration via Triangulation}\label{Integration Algorithm}
\endgroup
}
\end{algorithm}
\noindent Figure \ref{fig_triangulation} illustrates the partitioning of an arbitrary active subdomain into its triangular elements, facilitating the application of standard quadrature rules.

\subsubsection{Conforming Monte Carlo Integration}
While  Algorithm \ref{Integration Algorithm}  evaluates the spatial kernels via exact quadrature over the triangulated subdomains, an alternative approach is to apply Monte Carlo integration confined to the active subdomains $\Omega_{ij}$ and $\Omega_{i}$. We refer to this approach as Conforming Monte Carlo Integration.

Let $\phi(\bsx)$ denote a generic spatial kernel associated with the local operators $\mathcal{A}^{(l)}$ or $\mathcal{B}^{(l)}$ ($l=1, 2, 3$), and let $\psi(\bsx)$ be associated with $\mathcal{C}^{(l)}$ ($l=1, 2$). To construct the Conforming Monte Carlo estimators, we draw $N_{ij}$ and $N_i$ uniformly sampled evaluation points $\{\bsx_{ij,q}\}_{q=1}^{N_{ij}} \subset \Omega_{ij}$ and $\{\bsx_{i,q}\}_{q=1}^{N_i} \subset \Omega_i$, respectively. The corresponding domain integrals are approximated as follows
\begin{equation}
    \int_{\Omega_{ij}} \phi(\bsx) \, d\bsx \approx \frac{|\Omega_{ij}|}{N_{ij}} \sum_{q=1}^{N_{ij}} \phi(\bsx_{ij,q}), 
    \qquad
    \int_{\Omega_{i}} \psi(\bsx) \, d\bsx \approx  \frac{|\Omega_{i}|}{N_{i}} \sum_{q=1}^{N_{i}} \psi(\bsx_{i,q}),
\end{equation}
where $|\Omega_{ij}|$ and $|\Omega_{i}|$ denote the exact geometric areas (or volumes) of the clipped polygons. The key distinction from standard Monte Carlo integration is that the active subdomains $\Omega_{ij}$ and $\Omega_{i}$, together with their geometric measures, depend implicitly on the network parameters through the clipping hyperplanes. By retaining this dependence within the computational graph, automatic differentiation yields the correct gradients with respect to the network parameters. This approach incurs the standard stochastic convergence rate, i.e., $\mathcal{O}(1/\sqrt{N_{ij}})$, making it less accurate than the exact quadrature, but it provides a parallelizable and computationally straightforward alternative.

\subsection{Computational Cost of the Geometric Preprocessing}
Each clipping operation against a half-space costs $\mathcal{O}(v)$, where $v$ denotes the vertex count of the polygon being clipped, and adds at most one new vertex. The single-neuron region $\Omega_i$ is obtained by one clipping of $\Omega$, and the pairwise region $\Omega_{ij} = \Omega_i \cap \mbf{X}_j^+$ by one further clipping of $\Omega_i$. Hence, every region has at most $m+2$ vertices, where $m$ is the number of vertices of $\Omega$, and it is constructed at cost $\mathcal{O}(1)$. Since there are $n$ single-neuron regions and $\mathcal{O}(n^2)$ pairwise regions, the total preprocessing cost is $\mathcal{O}(n^2)$ per training step.

\section{Training Strategies}\label{sec:Training Strategies}

In this section, we investigate three distinct algorithmic strategies for 
training the neural network by minimizing the energy functional 
$\mathcal{J}(u_{\bstheta})$ defined in \eqref{energy_function}.

\subsection{Full-Batch Energy Minimization}
In the full-batch setting, the training process is formulated as the direct minimization of the global energy functional $\mathcal{J}(u_{\bstheta})$ over the entire computational domain. At each training epoch, the complete set of contributions to the energy functional is assembled using the matrix-vector formulation introduced in \eqref{J_matrix_fom1}. This results in a deterministic optimization procedure, where the gradient $\nabla_{\bstheta} \mathcal{J}(u_{\bstheta})$ reflects the full interaction between all neurons and all localized regions of the domain. While this approach ensures that the optimizer follows the exact gradient of the energy functional, it incurs an $\mathcal{O}(n^2)$ computational cost at each iteration, since all pairwise regions $\Omega_{ij}$ must be constructed and integrated to assemble the full stiffness matrix, making it computationally intensive for large-scale problems.
See Algorithm \ref{training_Algorithm_full_batch}.

\noindent
\begin{algorithm}[H]
\SetAlgoLined
\textbf{Initialization:} Initialize neural network parameters 
$\bstheta^{(0)} = (\bsc, W, \bsb)$ and set learning rate $\eta$;

\For{each epoch $k$}{
    Compute the energy functional 
    $\mathcal{J}(u_{\bstheta})$ using the matrix-vector form in \eqref{J_matrix_fom1};

    Compute the gradient with respect to the trainable parameters $$\mathbf{g} = \nabla_{\bstheta} \mathcal{J}(u_{\bstheta});$$

    Update the trainable parameters using the chosen optimizer.
}
\caption{Full-Batch Energy Minimization}
\label{training_Algorithm_full_batch}
\end{algorithm}

\subsection{Neuron-Wise Stochastic Energy Minimization}
To mitigate the computational overhead associated with global system assembly, we propose a second algorithm that decomposes the optimization process into localized energy contributions. Specifically, the matrix-vector formulation presented in Corollary \ref{corollary_matrix_form} can be reformulated as follows
\begin{equation}\label{J_matrix_fom_local_i}
\mathcal{J}(u_{\bstheta}) = \sum_{i=1}^n \mathcal{J}_i(u_{\bstheta}), \qquad \mathcal{J}_i(u_{\bstheta}) = \frac{1}{2} \bsr_i^\top \mathcal{K}_i \bsr - \bsr_i^\top \mathcal{F}_i,
\end{equation}
where $\bsr_i = [\bsz_i^\top, y_i]^\top \in \mathbb{R}^{d+1}$,  $\mathcal{K}_i$ denotes the $i$-th row block of the global stiffness matrix and $\mathcal{F}_i$ represents the local load contribution associated with the $i$-th neuron. To further enhance efficiency, optimization is performed using a neuron-wise update scheme. At each epoch $k$, we define a stochastic estimator of the energy functional $\mathcal{J}$ by selecting a uniformly random subset of indices $\mathbf{I}_k \subset \{1, \dots, n\}$ of cardinality $r\leq n$, where $r$ is a fixed constant. Following \eqref{J_matrix_fom_local_i}, the neuron-wise stochastic energy is defined as
\begin{equation}\label{J_matrix_fom_local_i_1}
\tilde{\mathcal{J}}_k(u_{\bstheta}) = \frac{n}{r}\sum_{i\in \mathbf{I}_k} \mathcal{J}_i(u_{\bstheta}).
\end{equation}

\noindent \begin{algorithm}[H]\SetAlgoLined\textbf{Initialization:} Initialize neural network parameters $\bstheta^{(0)} = (\bsc, W, \bsb)$, set $r$, and set learning rate;\\
\For {each epoch $k$}{
     Select a random subset $\mathbf{I}_k\subset \{1, \dots, n\}$ of size $r$;\\
{Compute the stochastic energy functional $\tilde{\mathcal{J}}_k$ using the matrix-vector form in \eqref{J_matrix_fom_local_i_1};\\
Compute the gradient with respect to the trainable parameters $${\mathbf{g}}_k = \nabla_{\bstheta} \tilde{\mathcal{J}}_k(u_{\bstheta});$$
\\
Update the trainable parameters vector using an optimizer.}}\caption{Neuron-Wise Stochastic Energy Minimization}\label{training_Algorithm_stochastic}\end{algorithm}

\noindent To provide theoretical justification for Algorithm \ref{training_Algorithm_stochastic}, we explain why updating the global parameter vector ${\bstheta}$ based on the stochastic energy functional $\tilde{\mathcal{J}}_k$ is mathematically sound. In the Deep Ritz framework, the solution $u_{\bstheta}$ is sought by minimizing the energy functional $\mathcal{J}(u)$, which reaches its minimum when the variation with respect to all basis functions vanishes. Each hidden neuron $i$ acts as a basis function supported on its active subdomain $\Omega_i$. Because the local contribution $\mathcal{J}_i$ explicitly includes the interaction with all neurons, the stochastic gradient $\nabla_{\bstheta} \mathcal{J}_i$ already contains information about the global state of the network. Updating the parameters based on a random subset of neurons allows the optimizer to perform localized corrections on a per-region basis, where the variance introduced by this stochastic sampling may act as a natural regularizer, potentially preventing the network from settling into poor local minima. This randomized strategy is justified by the linearity of the gradient operator. Since $\mathbf{I}_k$ is drawn uniformly at random among all subsets of $\{1,\dots,n\}$ of cardinality $r$, each index $i$ belongs to $\mathbf{I}_k$ with probability $r/n$, which does not depend on $\bstheta$. By linearity of expectation,
\begin{equation}
\mathbb{E}[\mathbf{g}_k]
    = \frac{n}{r}\sum_{i=1}^{n}\frac{r}{n}\,\nabla_{\bstheta}\mathcal{J}_i(u_{\bstheta})
    = \nabla_{\bstheta}\mathcal{J}(u_{\bstheta}).
\end{equation}
The expected stochastic gradient is therefore the true gradient, so the descent direction is preserved in expectation.

\subsection{Double Neuron-Wise Stochastic Energy Minimization}

The third algorithm extends Algorithm~\ref{training_Algorithm_stochastic}
by introducing a second level of stochasticity into the pairwise
interaction terms of each local energy contribution. Using the block
decomposition $\mathcal{K}_i\bsr = \sum_{j=1}^{n}\mathcal{K}_{ij}\bsr_j$,
the local energy \eqref{J_matrix_fom_local_i} can be written as
\begin{equation}\label{J_local_expanded}
    \mathcal{J}_i(u_{\bstheta})
    = \sum_{j=1}^{n}\frac{1}{2}\,\bsr_i^{\top}\mathcal{K}_{ij}\bsr_j
      - \bsr_i^{\top}\mathcal{F}_i,
\end{equation}
where $\mathcal{K}_{ij}$ denotes the $(i,j)$ block of the global stiffness matrix, encoding the pairwise interaction between neurons $i$ and $j$. 
At each epoch $k$, an outer uniformly random batch $\mathbf{I}_k \subset \{1,\dots,n\}$ of fixed size $r$ is selected as before. 
In addition, a single subset $\mathbf{J}_k$ of size $\bar r \le n$ is drawn uniformly at random, independently of $\mathbf{I}_k$, and shared by all $i \in \mathbf{I}_k$. Defining
\begin{equation}\label{J_matrix_fom_local_ij}
    \tilde{\mathcal{J}}_{ij}(u_{\bstheta})
    \;=\; \frac{1}{2}\,\bsr_i^{\top}
      \mathcal{K}_{ij}\bsr_j, \qquad \tilde{\mathcal{J}}_{i}(u_{\bstheta})
    \;=\; \frac{n}{\bar{r}}\Big(\sum_{j\in\mathbf{J}_{k}}\tilde{\mathcal{J}}_{ij}(u_{\bstheta})\Big)
      \;-\; \bsr_i^{\top}\mathcal{F}_i,
\end{equation}
the doubly stochastic energy of epoch $k$ is
\begin{equation}\label{J_matrix_fom_double_epoch}
\tilde{\tilde{\mathcal{J}}}_{k}(u_{\bstheta})
    \;= \frac{n}{r}\; \sum_{i\in\mathbf{I}_k}\tilde{\mathcal{J}}_{i}(u_{\bstheta}).
\end{equation}
Herein, note that the source term $\bsr_i^{\top}\mathcal{F}_i$ is retained in full, since it involves only local data associated with neuron $i$ and carries no dependence on the interaction index $j$. See Algorithm \ref{training_Algorithm_double_neuron_wise}.

\begin{algorithm}[H]
\SetAlgoLined
\textbf{Initialization:} Initialize neural network parameters
$\bstheta^{(0)} = (\bsc, W, \bsb)$ and set $r$, $\bar{r}$, and learning rate\;
\For{each epoch $k$}{
    Select a random subset $\mathbf{I}_k \subset \{1,\dots,n\}$ of
    size $r$\;
    Select a random subset $\mathbf{J}_{k} \subset \{1,\dots,n\}$
        of size $\bar{r}$\;
    \For{each neuron $i \in \mathbf{I}_k$}{
    \For{each neuron $j \in \mathbf{J}_k$}{
        Compute the pairwise term $\tilde{\mathcal{J}}_{ij}$
        via \eqref{J_matrix_fom_local_ij}\;
    }}
    Assemble $\tilde{\tilde{\mathcal{J}}}_{k}$
    via \eqref{J_matrix_fom_double_epoch}\;
    Compute the gradient with respect to the trainable parameters:
    \[
        \tilde{\mathbf{g}}_{k}
        = \nabla_{\bstheta}\,\tilde{\tilde{\mathcal{J}}}_{k}(u_{\bstheta});
    \]
    Update the trainable parameters $\bstheta$ using an optimizer\;
}
\caption{Double Neuron-Wise Stochastic Energy Minimization}
\label{training_Algorithm_double_neuron_wise}
\end{algorithm}

\noindent We now analyze the expected stochastic gradient of
Algorithm~\ref{training_Algorithm_double_neuron_wise}. Fix a neuron $i$. Each
index $j$ belongs to $\mathbf{J}_k$ with probability $\bar r/n$. By linearity
of expectation,
\begin{equation}
\mathbb{E}\bigl[\nabla_{\bstheta}\tilde{\mathcal{J}}_i\bigr]
  = \frac{n}{\bar r}\sum_{j=1}^{n}\frac{\bar r}{n}\,\nabla_{\bstheta}\tilde{\mathcal{J}}_{ij}
     - \nabla_{\bstheta}\bigl(\bsr_i^{\top}\mathcal{F}_i\bigr)
  = \sum_{j=1}^{n}\nabla_{\bstheta}\tilde{\mathcal{J}}_{ij}
     - \nabla_{\bstheta}\bigl(\bsr_i^{\top}\mathcal{F}_i\bigr)
  = \nabla_{\bstheta}\mathcal{J}_i(u_{\bstheta}).
\end{equation}
Each index $i$ belongs to $\mathbf{I}_k$ with probability $r/n$. Since
$\mathbf{I}_k$ and $\mathbf{J}_k$ are independent, we obtain
\begin{equation}
\mathbb{E}\bigl[\tilde{\mathbf g}_k\bigr]
  = \frac{n}{r}\sum_{i=1}^{n}\frac{r}{n}\,\mathbb{E}\bigl[\nabla_{\bstheta}\tilde{\mathcal{J}}_i\bigr]
  = \sum_{i=1}^{n}\nabla_{\bstheta}\mathcal{J}_i(u_{\bstheta})
  = \nabla_{\bstheta}\mathcal{J}(u_{\bstheta}).
\end{equation}
Hence, as in Algorithm~\ref{training_Algorithm_stochastic}, the expected
gradient equals the true gradient, and the descent direction is preserved in
expectation.

\vspace{0.5cm}

\begin{remark}

When a neuron's activation hyperplane leaves the computational domain $\Omega$, two situations arise. If $\Omega$ lies entirely on the inactive side, the neuron vanishes on $\Omega$. Its activation and all its parametric derivatives are zero, so it receives no gradient and remains inactive permanently. If, instead, $\Omega$ lies entirely on the active side, the neuron is affine on $\Omega$. It still receives gradient information, but it creates no activation boundary inside $\Omega$ and does not refine the partition of $\Omega$ into active subdomains. To mitigate these issues,
we propose augmenting the energy functional $\mathcal{J}(u_{\bstheta})$
with a geometric penalty grounded in the Hyperplane Separation
Theorem~\cite{vandenberghe2004convex}, which enforces that every
ReLU boundary intersects $\Omega$.

For a convex polygonal domain with vertex set $\mathcal{P}$, we introduce the normalized signed distance from each vertex $\bsp$ to the $j$-th hyperplane
\begin{equation}
    v_{j}(\bsp) = \frac{\bsw_j \cdot \bsp + b_j}{\|\bsw_j\|},
    \qquad j = 1, \ldots, n, \quad \bsp \in \mathcal{P}.
\end{equation}
The geometric penalty for the $j$-th neuron is then defined as
\begin{equation}\label{Reg}
    \mathcal{R}_j(\bstheta) =
    \Bigl[\mathrm{ReLU}\,\Bigl(\epsilon -
      \max_{\bsp \in \mathcal{P}} v_{j}(\bsp)\Bigr)\Bigr]^{2}
    +
    \Bigl[\mathrm{ReLU}\,\Bigl(\epsilon +
      \min_{\bsp \in \mathcal{P}} v_{j}(\bsp)\Bigr)\Bigr]^{2},
    \quad j = 1, \ldots, n,
\end{equation}
where $\epsilon \geq 0$ is a small safety margin. By construction,
$\mathcal{R}_j$ vanishes whenever the $j$-th ReLU boundary intersects $\Omega$ with margin $\epsilon$. If the boundary drifts entirely outside $\Omega$, $\mathcal{R}_j$ becomes strictly positive, generating a restoring gradient that pulls the boundary back into $\Omega$. The regularized training loss is
accordingly defined as
\begin{equation}\label{Loss with reg}
    \mathcal{L}_{\mathrm{reg}}(\bstheta)
    = \mathcal{J}(u_{\bstheta})
    + \frac{\lambda}{n}\sum_{j=1}^{n}\mathcal{R}_j(\bstheta),
\end{equation}
where $\lambda > 0$ is the regularization weight. 
Note that the regularized loss \eqref{Loss with reg} does not prevent the optimizer from discarding unnecessary neurons: it can drive their outer weights to zero, i.e., $c_j \to 0$, while their hyperplanes remain inside $\Omega$.
\end{remark}

\subsection{Computational costs}
\label{sec:cost}
This subsection evaluates the computational cost of the proposed method under each training strategy and standard Monte Carlo sampling.  For simplicity of presentation, we focus on the two-dimensional case in this study. However, this analysis can be straightforwardly extended to higher dimensions.

In the proposed framework, for each neuron $i$, we construct the active region $\Omega_i$, and for each pair $(i, j)$, we construct the pairwise intersection $\Omega_{ij} = \Omega_i \cap \Omega_j$. Given that the computational domain $\Omega$ is a convex polygon with a bounded number of vertices, both $\Omega_i$ and $\Omega_{ij}$ remain convex polygons with bounded numbers of vertices. Consequently, the geometric construction and evaluation of any single region $\Omega_i$ or $\Omega_{ij}$ incurs a constant cost of $\mathcal{O}(1)$. For notational convenience in subsequent complexity analyses, we define the following unit costs:
\begin{itemize}
  \item $c_{ij} = \mathcal{O}(1)$: the computational cost of evaluating the integral over one active region $\Omega_{ij}$;
  \item $c_{i} = \mathcal{O}(1)$: the computational cost of evaluating the integral over one active region $\Omega_i$;
  \item $c = \mathcal{O}(1)$: the baseline cost associated with computing the gradient with respect to one network parameter.
\end{itemize}

\paragraph{Standard Monte Carlo:} 
It samples directly on $\Omega$ without identifying active sub-polygons. Therefore, its geometric computational cost per iteration is 
\begin{equation}
T_{\rm MC} = 4n c = \mathcal{O}(n),
\end{equation}
because in the two-dimensional network \eqref{NN-definition}, we have $4n$ network parameters.

\paragraph{Full-batch:} According to the fact that $\Omega_{ij}=\Omega_{ji}$, a single
full-batch update uses all pairs $i<j$. The total geometric computational cost per epoch is thus given by
\begin{equation}
T_{\rm FB}=\frac{n(n-1)}{2}\,c_{ij} + n\,c_{i} + 4n\,c=\mathcal{O}(n^2)
\end{equation}

\paragraph{Neuron-Wise Stochastic Energy Minimization:} 
As outlined in Algorithm~\ref{training_Algorithm_stochastic}, the neuron-wise strategy selects at each iteration a subset of $r$ neurons, drawn uniformly at random. 
Therefore, the computational cost is given by
\begin{equation}
T_{\rm NW} \approx \frac{r(2n-r-1)}{2} \, c_{ij} + r \, c_{i} + 4n\, c = \mathcal{O}(n)
\end{equation}

\paragraph{Double Neuron-Wise Stochastic Energy Minimization:} 
As outlined in Algorithm~\ref{training_Algorithm_double_neuron_wise}, the double neuron-wise strategy draws two independent subsets of neurons uniformly at random of size $r$ and $\bar{r}$. 
Therefore, the computational cost per iteration is given by
\begin{equation}
T_{\rm DNW} \approx \left(r\bar{r} - \frac{r\bar{r}}{n}\right) c_{ij} + r \, c_{i} + 4(r + \bar r - \frac{r\bar{r}}{n})c 
= \mathcal{O}(1)
\end{equation}
where $\frac{r\bar{r}}{n}$ is the expectation of the overlap of two subsets. By choosing sub-linear sample sizes $r, \bar{r} \ll n$, this formulation significantly accelerates training while preserving robust stochastic gradient estimates.

\section{Numerical Experiments}\label{sec:Numerical experiments}
We present three numerical examples to assess the performance of the proposed methods. The first example considers a single-neuron benchmark designed to expose the limitations of standard Monte Carlo sampling and surrogate-gradient approaches; in particular, we demonstrate that both methods can fail to recover the true minimizer, whereas exact integration and Conforming Monte Carlo reliably converge to the correct solution. The second example employs a multi-neuron network with a polynomial right-hand side to compare the approximation accuracy of all four methods: standard Monte Carlo, surrogate gradient, exact integration, and Conforming Monte Carlo. This example is complemented by a comprehensive benchmark study that evaluates the trade-off between computational runtime and approximation accuracy across the training strategies introduced in Section~\ref{sec:Training Strategies}. The third example addresses a problem with a non-polynomial right-hand side function. Since exact integration is not applicable in this setting, all results are obtained with Conforming Monte Carlo and reported under each of the three training strategies.

Throughout all experiments, the regularized loss \eqref{Loss with reg} is employed with fixed parameters $\epsilon = 0.05$ and $\lambda = 0.1$.

\subsection{Deep Ritz limitations with Standard Monte Carlo and Surrogate Gradient: A One-Neuron Study}\label{A One-Neuron Study}

This example compares the proposed approaches on a two-dimensional test case
and examines whether each method attains the optimal solution.

\noindent Let $\Omega = (0,1)^2 \subset \mathbb{R}^2$. We seek a solution $u$
satisfying~\eqref{eq_PDE} with $a = 0$. The source term $f(x,y)$ is given by
\begin{equation}
    f(x,y) = 10(x + y) - 10.
\end{equation}
To approximate the solution, we employ a boundary-conforming NN architecture, as defined in~\eqref{NN-definition} and~\eqref{eq_DBC_NN}, with a single hidden neuron ($n = 1$) and the cutoff bubble function $g(x,y) = xy(x-1)(y-1)$. Due to the positive homogeneity of the ReLU activation, scaling $(\bsw_1, b_1)$ by a positive factor is equivalent to scaling $c_1$. Without loss of generality, we eliminate this redundancy by constraining $\bsw_1$ to unit length ($\|\bsw_1\|_2 = 1$) and parametrizing it via an orientation angle $\alpha \in [0, 2\pi)$ as $\bsw_1 = (\sin\alpha,\cos\alpha)^\top$. This parametrization allows $\alpha$ to control the orientation of the activation hyperplane, while $c_1$ governs the overall amplitude. Under this formulation, the trial approximation takes the form
\begin{equation}\label{eq_NN_1neuron}
u_{\bstheta}(x,y)
  = c_1 \cdot \operatorname{ReLU}\,\bigl(x\sin\alpha + y\cos\alpha + b_1\bigr)
    \cdot g(x,y).
\end{equation}
In this setting, analytically minimizing the energy functional $\mathcal{J}(u_{\bstheta})$, as defined in~\eqref{energy_function}, yields the optimal parameters $\alpha^* \approx 0.7854$, $c_1^* \approx 7.4323$, and $b_1^* \approx -0.5988$.

\noindent All simulations were initialized at $\alpha = \pi/2$, $c_1 = 1$, and $b_1 = -0.7$. The Adam optimizer was used with a learning rate exponentially decayed from $10^{-1}$ to $10^{-2}$. For the stochastic methods, $10^3 \times 10^3$ uniformly distributed collocation points were sampled for standard Monte Carlo, surrogate gradient, and Conforming Monte Carlo.

\begin{figure}[H]
    \centering
    \begin{subfigure}{\textwidth}
        \centering
        \includegraphics[width=0.9\linewidth]{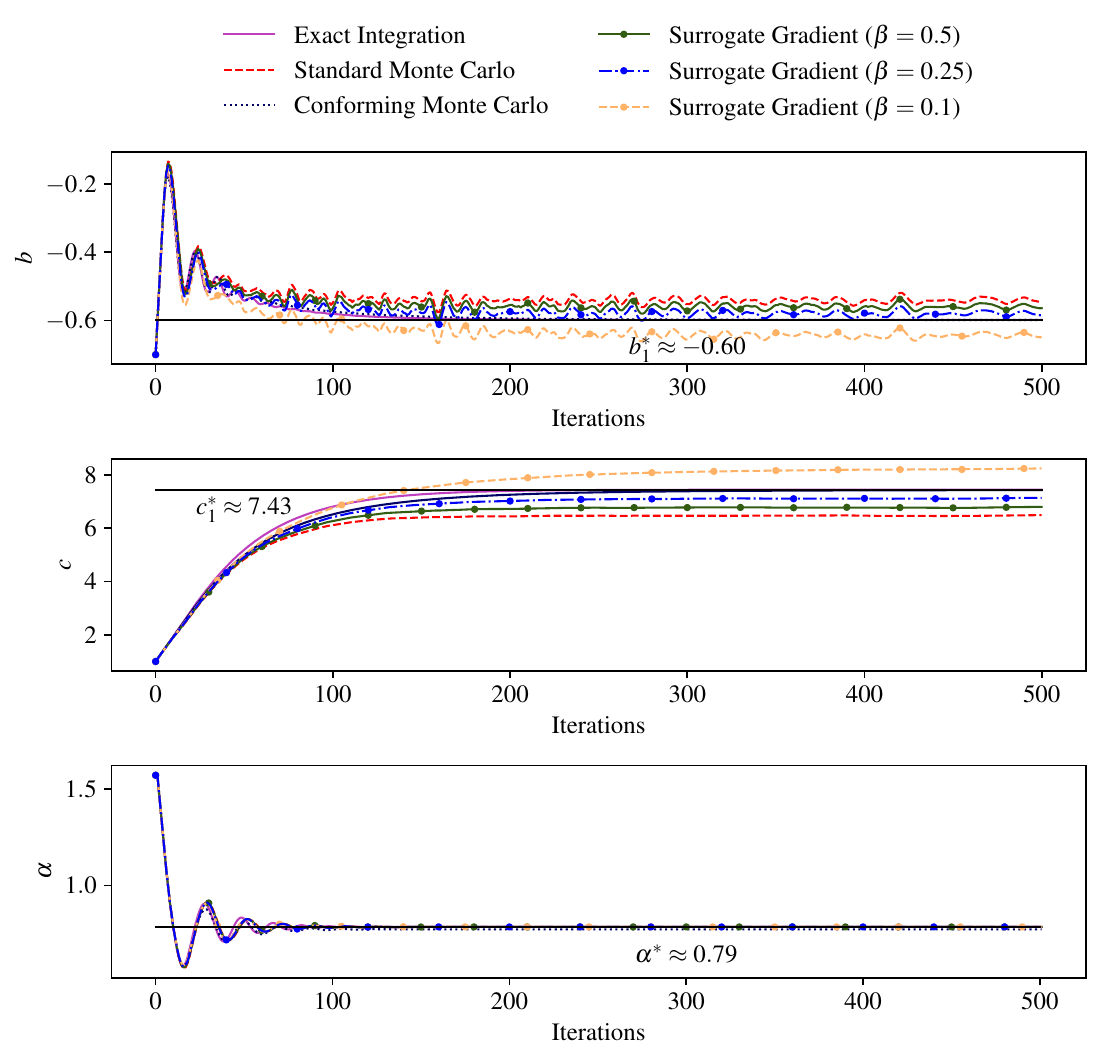}
        \caption{Optimization history over all $500$ iterations.}
        \label{fig:all_params_Example1}
    \end{subfigure}
    
    \vspace{1em} %
    
    \begin{subfigure}{\textwidth}
        \centering
        \includegraphics[width=1.0\linewidth]{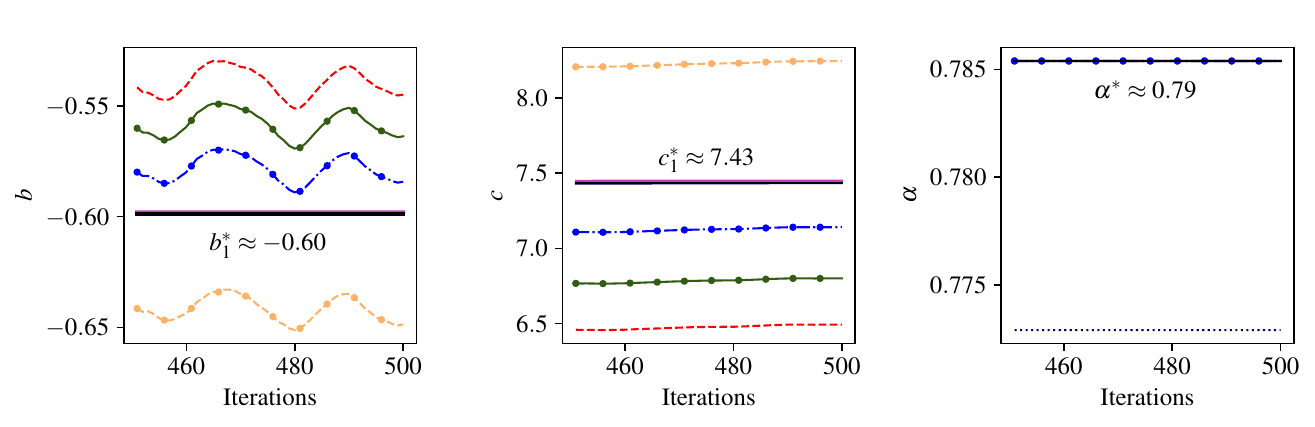}
        \caption{Magnified view of the optimization trajectory over the final $50$ iterations.}
        \label{fig:all_params_zoom_Example1}
    \end{subfigure}
    \caption{Convergence of the network parameters $\alpha$, $b_1$, and $c_1$ over $500$ training iterations for the single-neuron test case (Subsection~\ref{A One-Neuron Study}). Results are shown for exact integration, Conforming Monte Carlo, standard Monte Carlo, and the surrogate gradient method with smoothing parameters $\beta \in \{0.1, 0.25, 0.5\}$.}
    \label{fig:convergency_2D_Example1}
\end{figure}

\noindent Figure~\ref{fig:all_params_Example1}
compares the convergence histories of exact integration, Conforming Monte Carlo, standard Monte Carlo, and the surrogate gradient method under smoothing parameters $\beta \in \{0.1, 0.25, 0.5\}$ over $500$ training iterations. Figure \ref{fig:all_params_zoom_Example1} provides a zoomed view of the final $50$ iterations to highlight fine-scale differences in convergence behavior.

Although all methods rapidly identify the optimal orientation angle $\alpha^*$ within the first $100$ iterations, their behavior differs in learning the bias $b_1$ and the outer coefficient $c_1$. Both exact integration and Conforming Monte Carlo converge to the true analytical minimizers $b_1^*$ and $c_1^*$. In contrast, standard Monte Carlo and the surrogate gradient approaches fail to reach the optimal values. For the surrogate gradient method, the optimization trajectory is highly sensitive to the smoothing parameter $\beta$. Indeed, identifying a suitable value of $\beta$ is particularly challenging, as no universal criterion exists and $\beta$ must be calibrated empirically for each problem. Moreover, both Figures~\ref{fig:all_params_Example1} and \ref{fig:all_params_zoom_Example1} reveal significant instability in this approach, evidenced by persistent oscillations along the optimization trajectory. These results confirm that without explicitly resolving the activation boundary and its associated Dirac's delta sensitivities, standard sampling and smoothing strategies introduce biased gradients that prevent the optimizer from converging to the exact variational minimizer.

\subsection{A Multi-Neuron Benchmark}\label{A Multi-Neuron Benchmark}

This example pursues two objectives. First, we approximate the solution using a multi-neuron network and compare the four frameworks, exact integration, Conforming Monte Carlo, standard Monte Carlo, and the surrogate gradient method, in terms of approximation accuracy. We also examine the resulting spatial distribution of the ReLU activation hyperplane across $\Omega$ for each method. Second, we restrict attention to exact integration and Conforming Monte Carlo and evaluate their performance under the three training strategies introduced in Section~\ref{sec:Training Strategies}. For each strategy, we report approximation accuracy and computational runtime.

\noindent Consider the Poisson problem \eqref{eq_PDE} on the unit square domain $\Omega = (0, 1)^2 \subset \mathbb{R}^2$, subject to homogeneous Dirichlet boundary conditions, with the analytical solution given by
\begin{equation}
u(x, y) = (x^2 - x)(y^2 - y) xy.
\end{equation}
The corresponding source term, $f(x, y) = -\Delta u$, is derived analytically from the exact solution. We use a single-hidden-layer architecture as defined in \eqref{NN-definition} with 20 neurons ($n=20$). To strictly enforce the homogeneous Dirichlet boundary conditions, we employ the structure defined in \eqref{eq_DBC_NN} with $g(x, y) = (x^2 - x)(y^2 - y)$. Given the exact solution $u(x, y)$ and the boundary bubble function $g(x, y)$, the network $\mathcal{N}_{\boldsymbol{\theta}}(x, y)$ approximates the bilinear function $xy$. We employ the Adam optimizer with an exponential learning rate decay throughout training. The learning rate was scheduled to decay from $10^{-3}$ to $10^{-4}$ for the exact integration, $10^{-3}$ to $10^{-5}$ for the Conforming Monte Carlo, and from $10^{-1}$ to $10^{-2}$ for both the standard Monte Carlo and the surrogate gradient method. These ranges were chosen to ensure stable and efficient convergence behavior for each framework. For exact integration, we employ a $6 \times 6$ Gauss quadrature rule, whereas for the Monte Carlo approaches, we use a total of $10^3$ integration points, distributed uniformly at random over the entire domain $\Omega$.

\begin{figure}[htbp]
    \centering
    \includegraphics[width=0.8\textwidth]{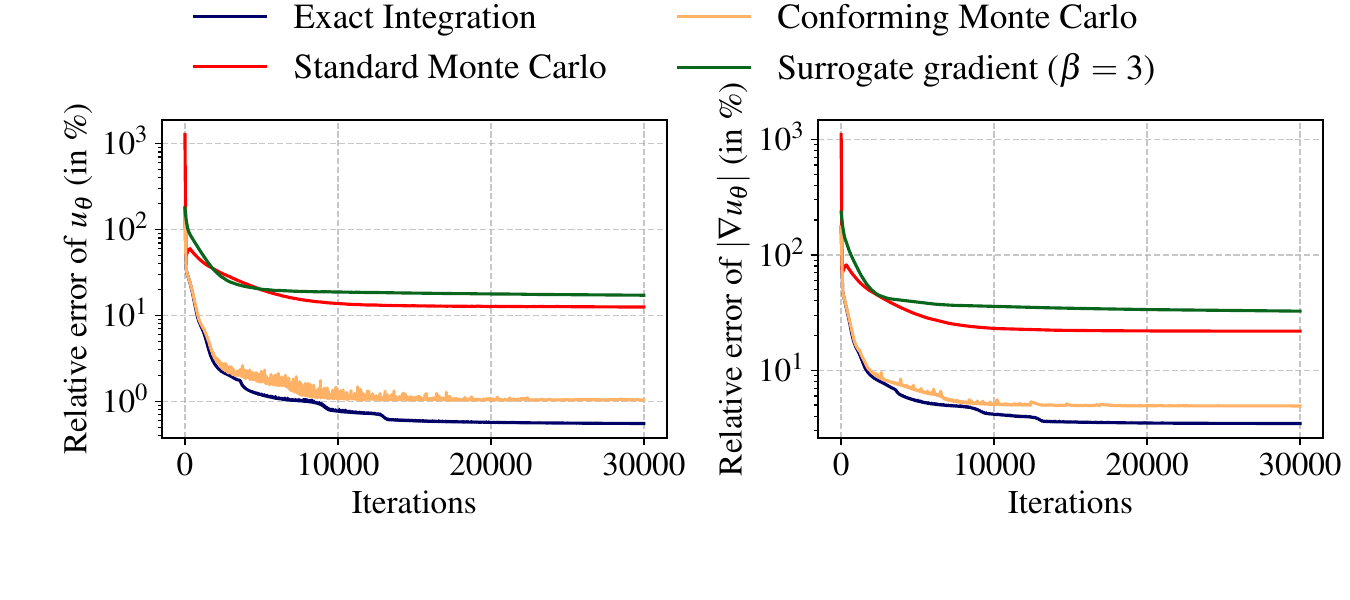}
    \vspace{-1cm}
    \caption{Evolution of relative $L^2$ errors (in $\%$) for the solution and its gradient during training for the benchmark in Subsection~\ref{A Multi-Neuron Benchmark}.}
    \label{fig_2D_errors_training_1}
\end{figure}

\paragraph{Convergence and activation behavior across integration methods:} Figure \ref{fig_2D_errors_training_1} illustrates the convergence histories of the relative $L^2$ error of the approximate solution $u_{\bstheta}$ and its gradient over $3\times10^4$ training iterations. Exact integration yields the most rapid and monotonic convergence, reducing the relative solution error to approximately $0.5\%$ and the gradient error to roughly $3.5\%$. Conforming Monte Carlo achieves comparable final accuracy; although stochastic oscillations are noticeable during intermediate iterations, the error decreases as the learning rate decays. In contrast, standard Monte Carlo and the surrogate gradient method remain at relative errors of about $12\%$ and $17\%$ for the solution $u_{\bstheta}$, and about $22\%$ and $33\%$ for its gradient $\nabla u_{\bstheta}$, respectively, after $3\times10^4$ iterations.

\begin{figure}[H]
    \centering
    \begin{subfigure}[b]{0.35\textwidth}
        \centering
    \includegraphics[width=\textwidth]{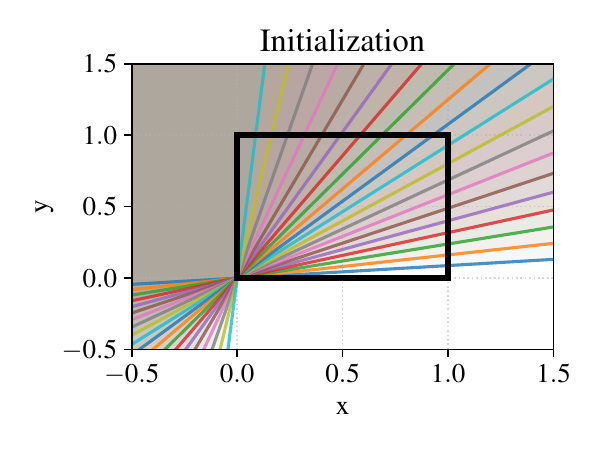}
        \vspace{-0.7cm}
    \end{subfigure}
    
    \begin{subfigure}[b]{0.35\textwidth}
        \centering
    \includegraphics[width=\textwidth]{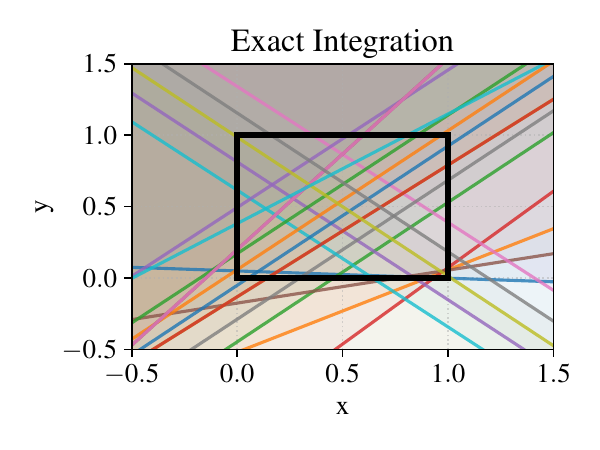}
        \vspace{-0.7cm}
    \end{subfigure}
    \begin{subfigure}[b]{0.35\textwidth}
        \centering
    \includegraphics[width=\textwidth]{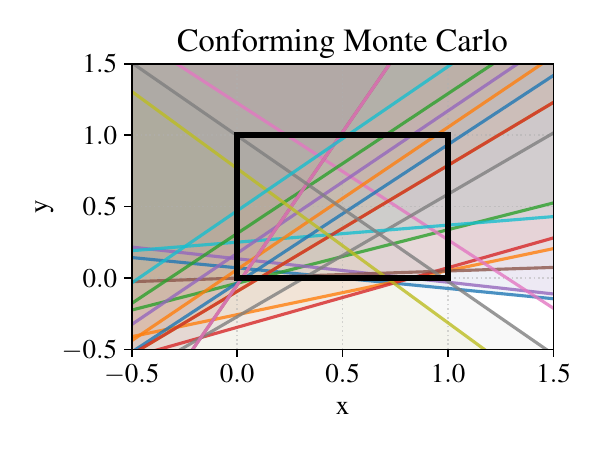}
        \vspace{-0.7cm}
        \label{fig:neuron_stochastic}
    \end{subfigure}
    \begin{subfigure}[b]{0.35\textwidth}
        \centering
    \includegraphics[width=\textwidth]{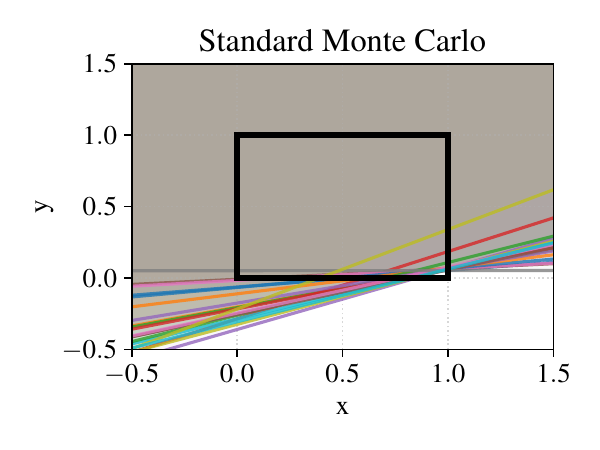}
        \vspace{-0.7cm}
        \label{fig:double_stochastic_1}
    \end{subfigure}
    \begin{subfigure}[b]{0.35\textwidth}
        \centering
        \includegraphics[width=\textwidth]{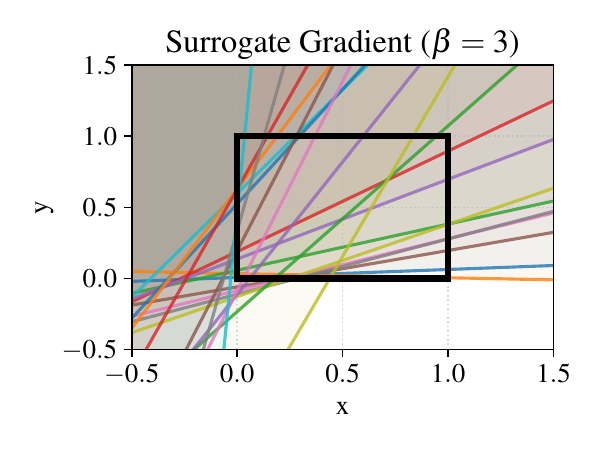}
        \vspace{-0.7cm}
        \label{fig:double_stochastic_2}
    \end{subfigure}

    \caption{ Initial and post-training ReLU activation patterns for the benchmark in Subsection \ref{A Multi-Neuron Benchmark}. Each colored line segment indicates the portion of a neuron's activation hyperplane within $\Omega$.}
    \label{fig_acrivate_ReLus}
\end{figure}

Figure \ref{fig_acrivate_ReLus} illustrates the ReLU activation supports within $\Omega$, comparing the initial configuration to the final distributions obtained by each method. The parameters were initialized such that the support of all twenty neurons intersects the domain at the start of training, providing a well-distributed initial basis for the optimization. As training progresses, the four methods exhibit substantially different behavior. Under standard Monte Carlo, nearly all breaking hyperplanes collapse toward the boundary of $\Omega$, as the method fails to capture gradient information from the ReLU breaking hyperplanes. The surrogate gradient method also performs poorly, and the breaking hyperplanes barely move from their initial positions. This indicates that the smoothed gradients lack the accuracy needed to drive effective parameter updates. In contrast, exact integration and Conforming Monte Carlo successfully distribute the breaking hyperplanes across the domain, allowing the optimizer to partition $\Omega$ into well-separated subdomains and achieve a more accurate approximation.

\begin{figure}[H]
    \centering
    \includegraphics[width=0.85\textwidth]{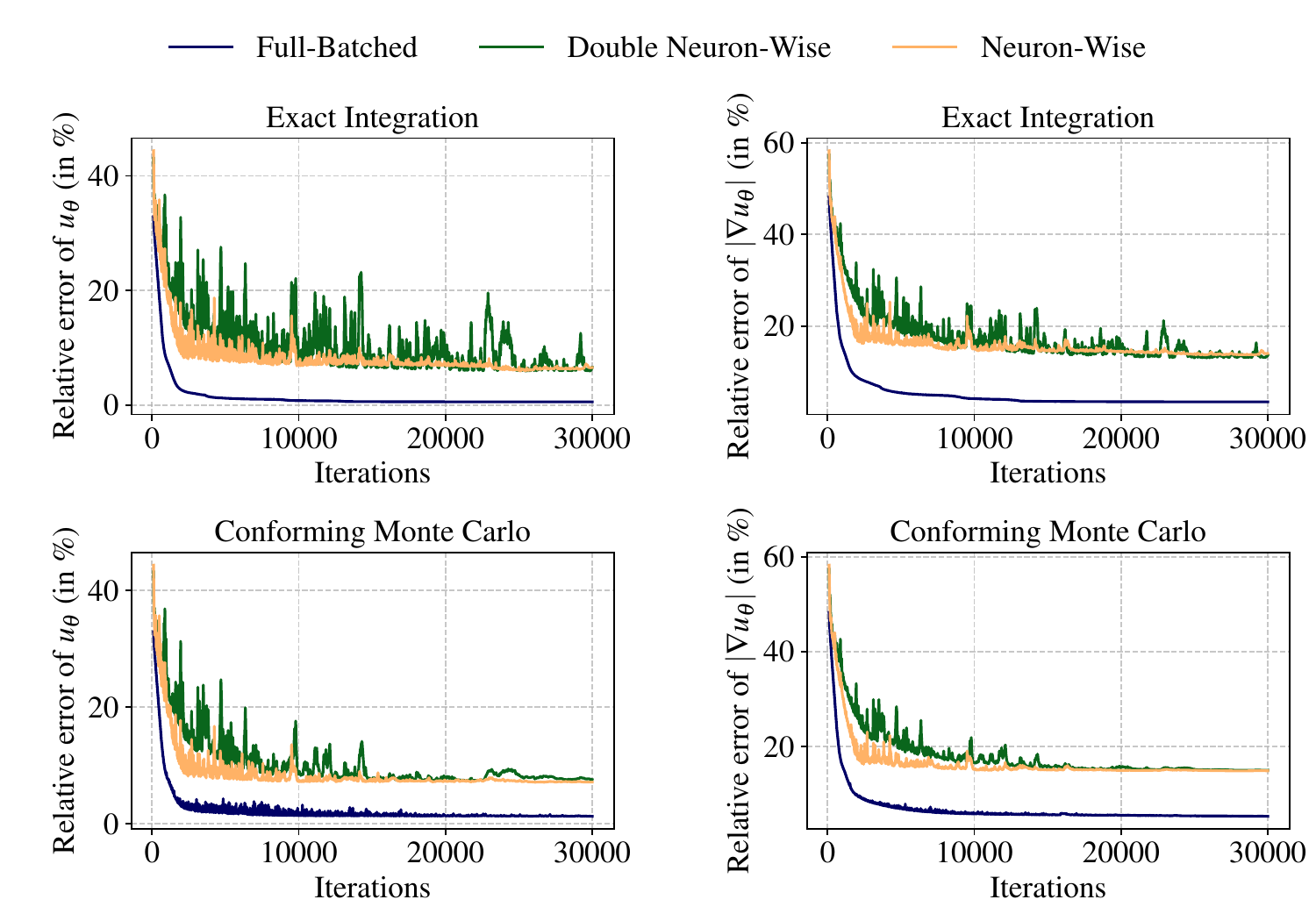}
    \caption{Evolution of relative $L^2$ errors (in $\%$) for the solution and its gradient during training for the benchmark in Subsection~\ref{A Multi-Neuron Benchmark}.}
    \label{fig_2D_errors_training}
\end{figure}

\paragraph{Comparison of Training Strategies:} Herein, we examine the training strategies described in
Section~\ref{sec:Training Strategies}. Figure~\ref{fig_2D_errors_training} illustrates the relative $L^2$ errors (in \%) for both the solution and its gradient, obtained using exact integration and the Conforming Monte Carlo framework under the three proposed training strategies: full-batch, neuron-wise, and double neuron-wise. For the neuron-wise strategy, we set the sample size $r = 10$, while for the double neuron-wise strategy, the sample sizes for both neuron subsets are selected as $r = \bar{r} = 10$. As shown, the full-batch strategy exhibits the most stable and monotonic convergence, rapidly reducing the error and achieving the lowest final error levels for both the solution and its gradient.
The neuron-wise approach follows a similar overall downward trajectory but introduces noticeable variance and minor oscillations during training; nevertheless, it converges to comparable final error levels. The double neuron-wise strategy, owing to its more heavily stochastic nature, exhibits more pronounced oscillations and the highest variance throughout training. Although its convergence is slower, it reduces the error over time and settles to an acceptable final error level. Overall, while stochastic strategies may exhibit slower convergence and a minor sacrifice in accuracy, they offer lower computational costs and may help escape poor local minima, establishing them as valuable alternatives to the full-batch approach. Furthermore, the relative performance of the three optimization strategies remains unchanged under both exact integration and the Conforming Monte Carlo framework. 

\begin{figure}[H]
    \centering
    \includegraphics[width=0.85\textwidth]{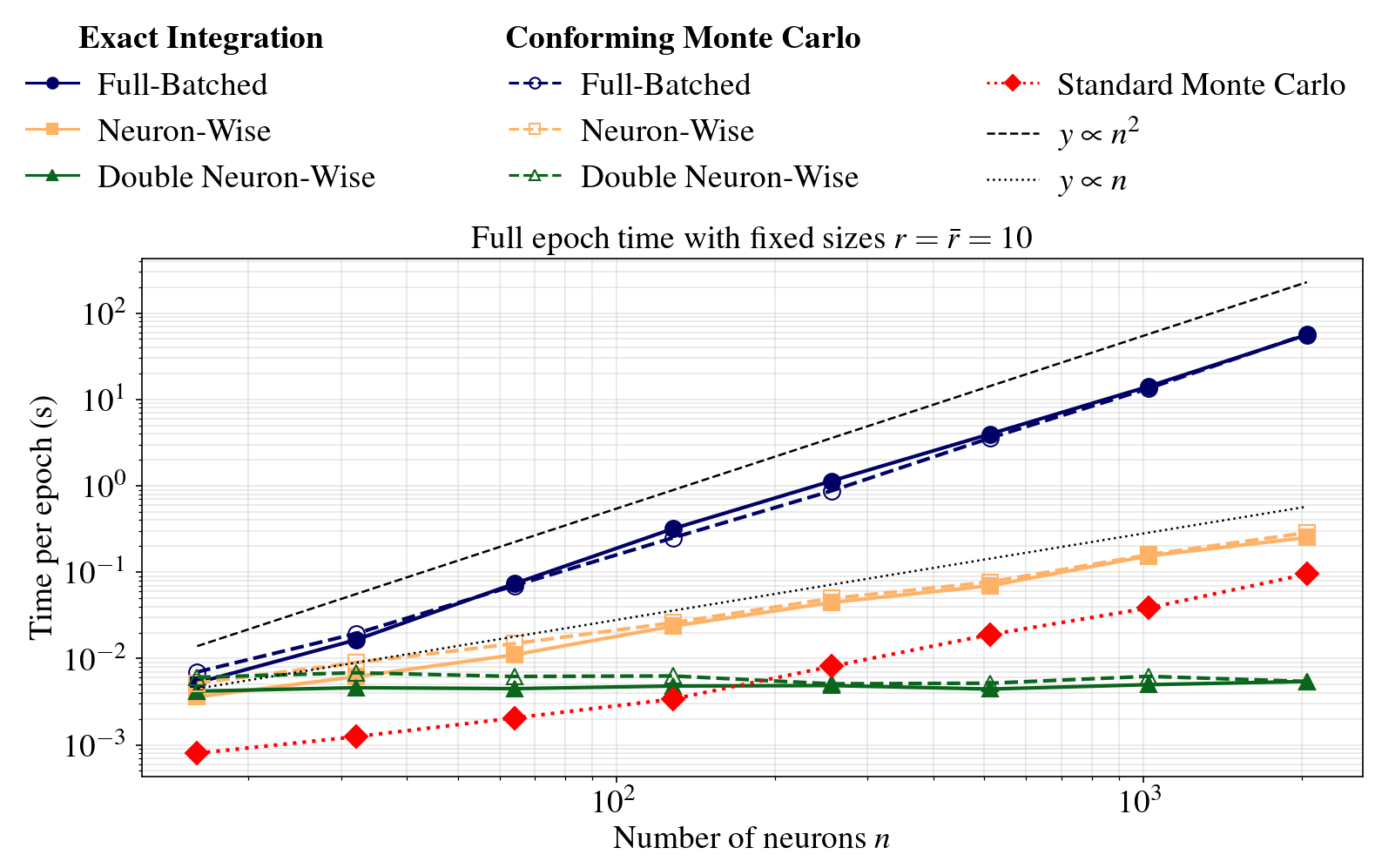}
    \caption{Time per epoch as a function of the number of neurons $n = 2^i$, $i = 4, \ldots, 11$, for exact integration, Conforming Monte Carlo, and standard Monte Carlo, using the full-batch, neuron-wise, and double neuron-wise training strategies (benchmark in Subsection~\ref{A Multi-Neuron Benchmark}). Both Monte Carlo methods use the same point density.}
    \label{fig_2D_time_reporting}
\end{figure}

\paragraph{Training Time:}
We now compare the computational time per epoch of the proposed training strategies. Figure~\ref{fig_2D_time_reporting} shows the time per epoch as a function of the number of neurons, $n=2^i$ for $i=4,\ldots,11$, for exact integration and Conforming Monte Carlo, each combined with the full-batch, neuron-wise and double neuron-wise strategies, and for standard Monte Carlo. All Monte Carlo based methods use the same integration point density of $10^4$ points per epoch. In this test, the subset sizes of the neuron-wise and double neuron-wise strategies are fixed at $r=\bar r=10$. The results agree with the theoretical analysis in Section~\ref{sec:cost}. The full-batch strategy grows quadratically with $n$ (order two), the neuron-wise strategy and standard Monte Carlo grow approximately linearly (order one), and the double neuron-wise strategy is essentially independent of $n$ (order zero). For each strategy, the curves of exact integration and Conforming Monte Carlo nearly coincide. For large $n$, the double neuron-wise strategy is the cheapest of all methods considered, and both stochastic strategies are substantially cheaper than the full-batch strategy. These results demonstrate the efficiency of the proposed training strategies.

\subsection{A Benchmark with Non-polynomial Right-hand Side}
\label{sec:non-polynomial}

\noindent We consider the two-dimensional Poisson equation~\eqref{eq_PDE} on
the unit square $\Omega = (0,1)^2 \subset \mathbb{R}^2$ with $a = 0$,
subject to homogeneous Dirichlet boundary conditions. The considered analytical solution is given by
\begin{equation}
    u(x, y) = \sin(\pi x)\sin(\pi y),
\end{equation}
with the corresponding source term $f(x, y) = -\Delta u$ that is derived analytically.

To strictly enforce homogeneous boundary conditions, we employ the structure~\eqref{eq_DBC_NN} with $g(x, y) = (x^2 - x)(y^2 - y)$. The trial solution is parametrized by an NN with a single hidden layer comprising $50$ neurons ($n = 50$), as defined in~\eqref{NN-definition}. Since the source term is non-polynomial, exact quadrature is not available; therefore, we use Conforming Monte Carlo integration to evaluate the loss function. For the neuron-wise strategy, we set the sample size $r = 20$, while for the double neuron-wise strategy, the sample sizes for both neuron subsets are selected as $r = \bar{r} = 20$. We also use a total of $10^3$ integration points, distributed uniformly at random over the entire domain $\Omega$. We train the network using the Adam optimizer with an exponential learning rate decay scheduled from $10^{-2}$ to $10^{-4}$.

\begin{figure}[htbp]
    \centering
    \includegraphics[width=1\textwidth]{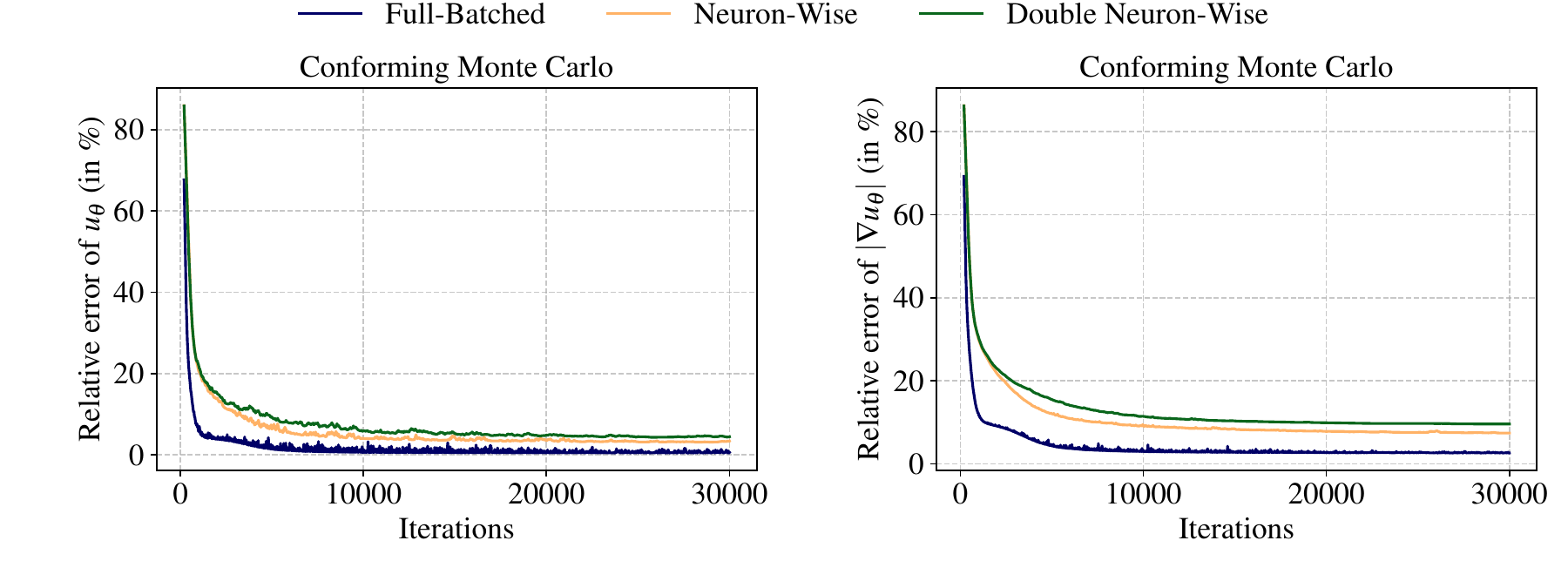}
    \caption{Evolution of relative $L^2$ errors for the solution and its
    gradient during training for the benchmark in
    Subsection~\ref{sec:non-polynomial}.}
    \label{fig_2D_errors_non-polynomial}
\end{figure}

Figure~\ref{fig_2D_errors_non-polynomial} presents the relative $L^2$ errors (in \%) for both the solution and its gradient over $3 \times 10^4$ training iterations, obtained using Conforming Monte Carlo under the three optimization strategies. All three strategies successfully reduce the error throughout training, confirming that Conforming Monte Carlo remains effective even when exact integration is unavailable. The full-batch strategy exhibits the most stable and monotonic convergence, achieving the lowest final error levels for both the solution and its gradient. The neuron-wise and double neuron-wise strategies introduce more pronounced oscillations during intermediate iterations, but converge to comparable final accuracy. These results demonstrate that the Conforming Monte Carlo framework provides a robust and flexible alternative to exact integration for problems with non-polynomial source terms, while preserving the accuracy advantages of correctly resolving the ReLU activation boundaries.

\section{Conclusion}\label{sec:Conclusion}

This work addresses a fundamental obstacle in training ReLU neural networks within the Deep Ritz Method or PINNs for first-order problems. Standard stochastic integration schemes, including Monte Carlo sampling and surrogate gradient methods, fail to capture the Dirac's delta contributions arising at the ReLU activation boundaries when differentiating the integrand of the loss. This failure produces incorrect gradient information regardless of the number of sample points, preventing reliable convergence to the correct variational minimizer.

To overcome this limitation, we developed an integration framework built around the Sutherland--Hodgman polygon clipping algorithm, which partitions the computational domain into convex active sub-polytopes where the network's spatial derivatives are piecewise constant. By explicitly resolving these sub-regions and their boundaries, the framework correctly captures the Dirac's delta distributions. The resulting integrals can be evaluated either to machine precision via triangular Gaussian quadrature or by a statistically consistent Conforming Monte Carlo approximation. In both cases, the loss is computed in matrix--vector form to preserve computational efficiency: all neuron interactions are computed at once with vectorized operations. It also splits the energy into neuron-wise blocks, which our stochastic training strategies rely on. Besides the Full-Batch strategy, we introduced two 
stochastic training strategies: Neuron-Wise Stochastic and 
Double Neuron-Wise Stochastic Energy Minimization. Each strategy offers a 
different trade-off between accuracy and computational cost. The stochastic 
strategies reduce the assembly cost per epoch relative to Full-Batch, while all 
three preserve the correct descent direction in expectation.

Three numerical experiments in two spatial dimensions validate the proposed framework. A one-neuron benchmark demonstrates that exact integration and Conforming Monte Carlo converge to the analytical minimizer, while standard Monte Carlo and surrogate gradient methods fail to do so. A multi-neuron benchmark confirms that the proposed frameworks accurately approximate both the solution and its gradient, while Monte Carlo fails once again. Finally, a third benchmark considers a problem whose source term is not a polynomial, making exact quadrature unavailable. In this setting, the Conforming Monte Carlo framework is applied under all three optimization strategies. In all cases, accurate and stable convergence is achieved, demonstrating the practical reach of the method beyond the polynomial setting.

Future work will address the extension of the framework to three-dimensional domains, deep ReLU networks, non-convex geometries, and more complex problems. In particular, an interesting future research work is applying our framework to PDE problems with sharp discontinuities or steep gradients, such as shock waves or interface problems. In such cases, one can design a network with ReLU activations only in the first hidden layer---to capture the discontinuous structure---and any smooth activation function in the remaining layers. Since the Dirac's delta distributions arise solely from the first ReLU layer, our framework correctly integrates them, ensuring reliable training.

\section*{Code Availability}
The source code for this paper is available at:

\url{https://github.com/Mathmode/ExactIntegration_ReLU-typeNN4PDEs}.

\section*{Acknowledgments}
SB, DP, and KvdZ are supported by the European Union’s Horizon Europe research and innovation programme under the Marie Skłodowska-Curie grant agreement No 101119556 (IN-DEEP). 

The research of IB and KvdZ was supported by the Engineering and Physical Sciences Research Council (EPSRC), UK, under Grant EP/W010011/1. 

Moreover, DP has received funding from the following Research Projects/Grants:
European Union’s Horizon Europe research and innovation programme
under the Marie Skłodowska-Curie Actions MSCA-SE-101342459 (MOST-PINN).
PID2023-146678OB-I00 funded by MICIU/AEI /10.13039/501100011033 and by
FEDER, EU;
BCAM Severo Ochoa accreditation of excellence CEX2021-001142-S funded
by MICIU/AEI/10.13039/501100011\\033;
Basque Government through the BERC 2022-2025 program;
EGOKIA (KK-2026/00050), funded by the Basque Government through ELKARTEK;
Consolidated Research Group MATHMODE (IT1866-26) of the UPV/EHU given
by the Department of Sciences, Universities, and Innovation of the
Basque Government;
BCAM-IKUR-UPV/EHU, funded by the Basque Government IKUR Strategy and
by the European Union NextGenerationEU/PRTR;

The authors thank Jamie M.~Taylor for early discussions regarding Monte Carlo integration for ReLU networks.

\section*{Appendix}
\appendix

\section{Proof of Proposition \ref{proposition_clipping}}\label{Proof of Proposition}

Let us consider the two-dimensional case, i.e., when
$\Omega \subset \mathbb{R}^2$ is a convex polygon. Note that $\Omega_i = \Omega \cap \mbf{X}_i^+$ is also convex, since it is the intersection of a half-space with a convex set. Throughout this section, we rely on the following standard characterization of vertices: a vertex in a convex polygon cannot be written as a convex
combination of two different points within the polygon; that is, $\bsp$ is a vertex of $\Omega_i$ if and only if there are no $\bsu, \bsv \in \Omega_i$, $\bsu \neq \bsv$, and no $\lambda \in (0,1)$ such that $\bsp = \lambda \bsu + (1-\lambda)\bsv$. Equivalently, every vertex is an extreme point of $\Omega_i$.

\medskip
\noindent First, let us consider the case where $\bsw_i \cdot \bsp + b_i > 0$ and $\bsp \in \mathcal{P}$. Since $\bsp$ is already a vertex of $\Omega$, it cannot be written as a convex combination of two distinct points of $\Omega$, and therefore it cannot be written as a convex combination of two distinct points in the smaller set $\Omega_i \subseteq \Omega$. Hence $\bsp$ is a vertex of $\Omega_i$.

\medskip
\noindent Now, let us consider the case where $\bsw_i \cdot \bsp + b_i = 0$ and
$\bsp \in \operatorname{Conv}(\bsp_k, \bsp_{k+1})$ for some $k$. Let 
$S_i:= \Gamma_i\cap \Omega$. Since $\Gamma_i$ contains a point of $\operatorname{int}(\Omega)$ and $\Omega$ is compact and convex, $S_i$ is a nondegenerate closed line segment whose interior lies in $\operatorname{int}(\Omega)$. In particular, the only points of $S_i$ on $\partial\Omega$ are its two endpoints. Since $\bsp \in \operatorname{Conv}(\bsp_k, \bsp_{k+1}) \subset \partial\Omega$ and $\bsw_i \cdot \bsp + b_i = 0$, we have $\bsp \in \Gamma_i \cap \partial\Omega$, so $\bsp$ is an endpoint of $S_i$.

Suppose that $\bsp$ can be written as 
\begin{equation}\label{p_def}
\bsp = \lambda \bsu + (1-\lambda)\bsv
\end{equation}
for some distinct points $\bsu, \bsv \in \Omega_i$ and $\lambda \in (0,1)$. Since $\bsu,\bsv \in \Omega_i$, we have that $\bsw_i \cdot \bsu + b_i \geq 0$ and $\bsw_i \cdot \bsv + b_i \geq 0$. Then, 
\begin{equation}\label{convex_proof}
0 = \bsw_i \cdot \bsp + b_i = \lambda\,(\bsw_i \cdot \bsu + b_i)
\;+\;
(1-\lambda)\,(\bsw_i \cdot \bsv + b_i),
\end{equation}
which implies $\bsw_i \cdot \bsu + b_i = \bsw_i \cdot \bsv + b_i = 0$, and hence $\bsu, \bsv \in S_i$. Since $\bsp$ is an endpoint of $S_i$ and $\bsu, \bsv \in S_i$, the endpoint $\bsp$ cannot be an interior convex combination of two points in $S_i$ unless both coincide with $\bsp$; therefore $\bsu = \bsv = \bsp$, contradicting $\bsu \neq \bsv$. Consequently, $\bsp$ cannot be expressed as a convex combination of two distinct points in $\Omega_i$, and hence $\bsp$ is a vertex of $\Omega_i$.

\small
\bibliographystyle{abbrv}
\bibliography{references}
\end{document}

%% file: StandardSetOfCommands.tex
\newcommand{\ReLU}{\operatorname{ReLU}}

\newlength{\bigfboxsep}
\newcommand{\bs}[1]{\boldsymbol{#1}} 

\newcommand{\bsb}{\bs{b}}
\newcommand{\bsc}{\bs{c}}

\newcommand{\bsp}{\bs{p}}

\newcommand{\bsr}{\bs{r}}

\newcommand{\bsu}{\bs{u}}
\newcommand{\bsv}{\bs{v}}
\newcommand{\bsw}{\bs{w}}
\newcommand{\bsx}{\bs{x}}
\newcommand{\bsy}{\bs{y}}
\newcommand{\bsz}{\bs{z}}

\newcommand{\bstheta}{\bs{\theta}}

\newcommand{\bsTheta}{\bs{\Theta}}

\newcommand{\mbf}[1]{\mathbf{#1}}

%% file: Figure_Base.tex
\begin{tikzpicture}[
punkt/.style={
   rectangle,
   rounded corners,
   draw=black, thick,
   text width=2em,
   minimum height=1em,
   text centeblack}
]

\Vertex[x=0,y=1.5,label=$ \bsx$,color=white,size=1,fontsize=\normalsize]{X}
\draw (0,0.4) node {Input};

\Vertex[x=2,y=3.5,label=ReLU$(\cdot)$,color=blue!20,size=1.2,fontsize=\scriptsize]{A11}
\Vertex[x=2,y=2,label=ReLU$(\cdot)$,color=blue!20,size=1.2,fontsize=\scriptsize]{A12}
\Vertex[x=2,y=-0.5,label=ReLU$(\cdot)$,color=blue!20,size=1.2,fontsize=\scriptsize]{A13}

\Edge[color=black, Direct,lw=0.5pt](X)(A11)
\Edge[color=black, Direct,lw=0.5pt](X)(A12)
\Edge[color=black, Direct,lw=0.5pt](X)(A13)

\Vertex[x=4,y=1.5,label=$\mathcal{N}_\theta (\bsx)$,color=white,size=1.2,fontsize=\normalsize]{B}
\draw (4,0.4) node {Output};

\Edge[color=black, Direct,lw=0.5pt](A11)(B)
\Edge[color=black, Direct,lw=0.5pt](A12)(B)
\Edge[color=black, Direct,lw=0.5pt](A13)(B)

\draw (5.5,2) node{ $ \downarrow $};
\draw (5.5, 2.7) node{Apply $g(\bsx)$};

\Vertex[x=7,y=1.5,label=$u_\theta(x)$,color=red!20,size=1.2,fontsize=\normalsize]{C}
\Edge[color=black, Direct,lw=0.5pt](B)(C)
\draw (7,0.4) node {Approximate};
\draw (7,0) node {solution};
\end{tikzpicture}

%% file: Clipping_polygons.tex
\begin{tikzpicture}[
punkt/.style={
   rectangle,
   rounded corners,
   draw=black, thick,
   text width=2em,
   minimum height=1em,
   text centered}
]
\coordinate (A) at (2,-0.8);
\coordinate (B) at (4,0);
\coordinate (C) at (3.5,3);
\coordinate (D) at (1,3.5);
\coordinate (E) at (-1,1.5);
\coordinate (F) at (0,0.74);
\coordinate (G) at (3.61,2.31);
\coordinate (H) at (1.72,2.58);
\draw[thick] (A) -- (B) -- (C) -- (D) -- (E) -- cycle;
\fill[pattern=north east lines, pattern color=gray!80] (A) -- (B) -- (G) -- (H) -- (F) -- cycle;
\draw[thick, red] (-1.4,-0.8) -- (3,4);
\node[red] at (4,3.5) {$\bsw_i \cdot \bsx + b_i \geq 0$};
\draw[thick, blue] (-1,3) -- (5,2.1);
\node[blue] at (5.2,1.6) {$\bsw_j \cdot \bsx + b_j \geq 0$};
    \coordinate (P1) at (-0.7,-0.03);
    \draw[->,thick, red] (P1) -- ++(0.45,-0.4);
    \coordinate (P2) at (-0.8,2.97);
    \draw[->,thick, blue] (P2) -- ++(-0.12,-0.55);    %
\foreach \point in {F, A, B, G, H}
    \fill[black] (\point) circle (2pt);
\node at (0.8,-0.5) {\large $\Omega$};

\end{tikzpicture}

%% file: Fig_triangles.tex
\begin{tikzpicture}[font=\small, scale=1.5]
    \coordinate (V1) at (0.5, 0);
    \coordinate (V2) at (2, -0.3);
    \coordinate (V3) at (2.3, 2);
    \coordinate (V4) at (1.2, 2.2);
    \coordinate (V5) at (0, 1.2);
    \coordinate (C) at (1.24, 1.04); %

    \draw[fill=blue!10, draw=blue!60, thick] (V1) -- (V2) -- (V3) -- (V4) -- (V5) -- cycle;
    
    \foreach \n in {1,2,3,4,5} {
        \draw[dashed, gray] (C) -- (V\n);
        \filldraw[black] (V\n) circle (1pt) node[anchor=south east] {$\bsp_{\n}$};
    }
    
    \filldraw[red] (C) circle (1pt) node[anchor=north] {$\bsx_c$};
    
    \draw[fill=red!20, opacity=0.5] (C) -- (V2) -- (V3) -- cycle;
    \node at (1.9, 1.1) {$\mathcal{T}_2$};
    \node at (1.2, 0.3) {$\mathcal{T}_1$};
    \node at (1.4, 1.5) {$\mathcal{T}_3$};
    \node at (0.3, 1.9) {$\Omega_{ij}$};
\end{tikzpicture}